# INTEGRATION BY PARTS FORMULA FOR LOCALLY SMOOTH LAWS AND APPLICATIONS TO SENSITIVITY COMPUTATIONS

By Vlad Bally, Marie-Pierre Bavouzet and Marouen Messaoud

*Université de Marne la vallée, INRIA Rocquencourt and IXIS and INRIA Rocquencourt*

We consider random variables of the form $F=f(V_1,\dots,V_n)$, where $f$ is a smooth function and $V_i, i\in\mathbb{N}$, are random variables with absolutely continuous law $p_i(y)\,dy$. We assume that $p_i$, $i=1,\dots,n$, are piecewise differentiable and we develop a differential calculus of Malliavin type based on $\partial \ln p_i$. This allows us to establish an integration by parts formula $E(\partial_i\phi(F)G)=E(\phi(F)H_i(F,G))$, where $H_i(F,G)$ is a random variable constructed using the differential operators acting on $F$ and $G$. We use this formula in order to give numerical algorithms for sensitivity computations in a model driven by a Lévy process.

**1. Introduction.** In recent years, following the pioneering papers [12, 13], much work concerning numerical applications of stochastic variational calculus (Malliavin calculus) has been carried out. This mainly concerns applications in mathematical finance: computation of conditional expectations (which appear in, e.g., American option pricing) and of sensitivities (the so-called *Greeks*). The models at hand are usually log-normal type diffusions and one may then use standard Malliavin calculus. Currently, there is increasing interest in jump-type diffusions (see, e.g., [7]) and one must then use the stochastic variational calculus corresponding to Poisson point processes. Such a calculus has already been developed (in [4] and [15]) concerning the noise coming from the amplitudes of the jumps and (in [6, 9, 18, 19, 21, 22]) concerning times. Recently, Bouleau (see [5]) established the so-called *error calculus* based on the Dirichlet forms language and showed that the approaches in [4] and [6] fit into this framework. Another point of view, based on chaos decomposition, may be found in [3, 10, 16, 17, 23].

Let us finally mention that the models considered in mathematical finance (e.g., Merton's model) may have both a diffusion component (driven by a Brownian motion) and a jump part (driven by a compound Poisson process). In this case, one may use the standard Malliavin calculus with respect to the Brownian motion increments after conditioning in a clear way with respect to the Poisson component. This is done in [8, 11, 20].

The aim of this paper is to give a concrete application of the Malliavin calculus approach to sensitivity computations (Greeks) for pure jump diffusion models. We give three examples: in the first, we use the Malliavin calculus with respect to the jump amplitudes and in the second, we differentiate with respect to the jump times. In the third, we differentiate with respect to both of them.

The basic tool is an integration by parts formula which is analogous to the one in the standard Malliavin calculus on the Wiener space. Here, we give an abstract approach which, in particular, permits us to treat in an unified way the derivatives with respect to the times and the amplitudes of the jumps of Lévy processes. More precisely, we consider functionals of a finite number of random variables $V_i, i = 1, \dots, n$. The only assumption is that for each $i = 1, \dots, n$, the conditional law of $V_i$ (with respect to $V_j$, $j \neq i$) is absolutely continuous with respect to the Lebesgue measure and the conditional density $p_i = p_i(\omega, y)$ is piecewise differentiable. Using integration by parts, one may settle the duality relation which represents the starting point in Malliavin calculus. However, some border terms will appear corresponding to the points at which $p_i$ is not continuous: for example, if $V_i$ has a uniform conditional law on $[0,1]$, the density is $p_i(\omega, y) = 1_{[0,1]}(y)$ and integration by parts produces border terms in 0 and in 1. There is a simple idea which permits us to cancel the border terms: we introduce in our calculus some weights $\pi_i$ which are null at the points of singularity of $p_i$—in the previous example, we may take $\pi_i(y) = y^\alpha(1-y)^\alpha$ with some $\alpha \in (0,1)$. We then obtain a standard duality relation between the Malliavin derivative and the Skorohod integral and the machinery established in the Malliavin calculus produces an integration by parts formula. But there is a difficulty hidden in this process: the differential operators involve the weights $\pi_i$ and their derivatives. In the previous example, $\pi_i'(\omega, y) = \alpha(y^{\alpha-1}(1-y)^\alpha - y^\alpha(1-y)^{\alpha-1})$. These derivatives blow up in the neighborhood of the singularity points and this produces some nontrivial integrability problems. We must therefore search an equilibrium between the speed of convergence to zero and the speed with which the derivatives of the weights blow up in the singularity points. This leads to a nondegeneracy condition which involves the weights and their derivatives.

The integration by parts formula is established in Section 2. Since numerical algorithms involve only functions of a finite number of variables, we do not develop here an infinite-dimensional Malliavin calculus, but restrict ourselves to simple functionals. In Section 3, we use the integration by parts formula in order to compute the Delta (derivative with respect to the initial condition) for European options based on an asset which follows a pure jump diffusion equation and in Section 4, we give numerical results.

## 2. Malliavin calculus for simple functionals.

2.1. *The frame.* We consider a probability space $(\Omega, \mathcal{F}, P)$, a sub-$\sigma$-algebra $\mathcal{G} \subseteq \mathcal{F}$ and a sequence of random variables $V_i, i \in \mathbb{N}$. We define $\mathcal{G}_i = \mathcal{G} \vee \sigma(V_j, j \neq i)$. Our aim is to establish an integration by parts formula for functionals of $V_i$, $i \in \mathbb{N}$, which is analogous to the one in Malliavin calculus. The $\sigma$-algebra $\mathcal{G}$ appears in order to describe all the randomness which is not involved in the differential calculus.

We will work on some set $A \in \mathcal{G}$ which will be fixed throughout this section. We denote by $L_{(\infty)}(A)$ the space of the random variables such that $\mathrm{E}(|F|^p \mathbf{1}_A) < \infty$ for all $p \in \mathbb{N}$, and $L_{(p+)}(A)$ will be the space of the random variables $F$ for which there exists some $\delta > 0$ such that $\mathrm{E}(|F|^{p+\delta} \mathbf{1}_A) < \infty$. We assume the following.

Hypothesis 2.1. $V_i \in L_{(\infty)}(A)$, $i \in \mathbb{N}$.

For each $i \in \mathbb{N}$, we consider some $k_i \in \mathbb{N}$ and some $\mathcal{G}_i$-measurable random variables

$$a_i(\omega) = t_i^0(\omega) < t_i^1(\omega) < \cdots < t_i^{k_i}(\omega) < t_i^{k_i+1}(\omega) = b_i(\omega).$$

We define

$$B_i(\omega) = \bigcup_{j=0}^{k_i} (t_i^j(\omega), t_i^{j+1}(\omega)).$$

Notice that we may take

$$a_i = -\infty \quad \text{and} \quad b_i = \infty.$$

We will work with functions defined on $(a_i(\omega), b_i(\omega))$ which are smooth except for the points $t_i^j, j = 1, \ldots, k_i$. We define $\mathcal{C}_k(B_i)$ to be the set of measurable functions $f : \Omega \times \mathbb{R} \to \mathbb{R}$ such that for every $\omega$, $y \to f(\omega, y)$ is $k$-times differentiable on $B_i(\omega)$ and for each $j = 1, \ldots, k_i$, the left-hand side and the right-hand side limits $f(\omega, t_i^j(\omega)-), f(\omega, t_i^j(\omega)+)$ exist and are finite [for $j = 0$ (resp., $j = k_i + 1$), we assume that the right-hand side (resp., the left-hand side) limit exists and is finite]. We define

$$\begin{aligned} \Gamma_i(f) = {} & \sum_{j=1}^{k_i} (f(\omega, t_i^j(\omega)-) - f(\omega, t_i^j(\omega)+)) \\ & + f(\omega, b_i(\omega)-) - f(\omega, a_i(\omega)+). \end{aligned} \tag{2.1}$$

For $f,g\in\mathcal{C}_1(B_i)$, the integration by parts formula gives

$$\int_{(a_i,b_i)} fg'(\omega,y)\,dy=\Gamma_i(fg)-\int_{(a_i,b_i)} f'g(\omega,y)\,dy, \tag{2.2}$$

so $\Gamma_i$ represents the contribution of the border terms—or, in other words, of the singularities of $f$ or $g$.

Let $n,k\in\mathbb{N}$. We denote by $\mathcal{C}_{n,k}$ the class of the $\mathcal{G}\times\mathcal{B}(\mathbb{R}^n)$-measurable functions $f:\Omega\times\mathbb{R}^n\to\mathbb{R}$ such that $I_i(f)\in\mathcal{C}_k(B_i)$, $i=1,\dots,n$, where

$$I_i(f)(\omega,y):=f(\omega,V_1,\dots,V_{i-1},y,V_{i+1},\dots,V_n).$$

For a multi-index $\alpha=(\alpha_1,\dots,\alpha_k)\in\{1,\dots,n\}^k$, we define

$$\partial_\alpha^k f=\frac{\partial^k}{\partial x_{\alpha_1}\cdots\partial x_{\alpha_k}}f.$$

Moreover, we denote by $\mathcal{C}_{n,k}(A)$ the space of functions $f\in\mathcal{C}_{n,k}$ such that for every $0\le p\le k$ and every $\alpha=(\alpha_1,\dots,\alpha_p)\in\{1,\dots,n\}^p$, $\partial_\alpha^p f(V_1,\dots,V_n)\in L_{(\infty)}(A)$.

The points $t_i^j, j=1,\dots,k_i$, represent singularity points of the functions at hand (note that $f$ may be discontinuous at $t_i^j$) and our main purpose is to establish a calculus adapted to such functions.

Our basic hypothesis is the following.

HYPOTHESIS 2.2. For every $i\in\mathbb{N}$, the conditional law of $V_i$ with respect to $\mathcal{G}_i$ is absolutely continuous on $(a_i,b_i)$ with respect to the Lebesgue measure. This means that there exists a $\mathcal{G}_i\times\mathcal{B}(\mathbb{R})$-measurable function $p_i=p_i(\omega,x)$ such that

$$\mathrm{E}(\Theta\psi(V_i)\mathbf{1}_{(a_i,b_i)}(V_i))=\mathrm{E}\biggl(\Theta\int_{\mathbb{R}}\psi(x)p_i(\omega,x)\mathbf{1}_{(a_i,b_i)}(x)\,dx\biggr)$$

for every positive, $\mathcal{G}_i$-measurable random variable $\Theta$ and every positive, measurable function $\psi:\mathbb{R}\to\mathbb{R}$.

We assume that $p_i\in\mathcal{C}_1(B_i)$ and $\partial_y\ln p_i(\omega,y)\in L_{(\infty)}(A)$.

In concrete problems, we consider random variables $V_i$ with conditional densities $p_i$ and we then take $t_i^j$, $i=0,\dots,k_{i+1}$, to be the singularities of $p_i$. This means that we choose $B_i$ in such a way that $p_i$ satisfies Hypothesis 2.2 on $B_i$. This is the significance of $B_i$ (in the case where $p_i$ is smooth on the whole $\mathbb{R}$, we may choose $B_i=\mathbb{R}$).

For each $i\in\mathbb{N}$, we consider a $\mathcal{G}_i\times\mathcal{B}(\mathbb{R})$-measurable and positive function $\pi_i:\Omega\times\mathbb{R}\to\mathbb{R}_+$ such that $\pi_i(\omega,y)=0$ for $y\notin(a_i,b_i)$ and $\pi_i\in\mathcal{C}_1(B_i)$. We assume the following.

HYPOTHESIS 2.3. $\pi_i\in L_{(\infty)}(A)$ and $\pi_i'\in L_{(1+)}(A)$.

These will be the weights used in our calculus. In the standard Malliavin calculus, they appear as renormalization constants. On the other hand, $p_i$ may have discontinuities at $t_i^j, j = 1, \ldots, k_i$, and this will produce some border terms in the integration by parts formula; see (2.2). We may choose $\pi_i$ in order to cancel these border terms (as well as the border terms in $a_i$ and $b_i$).

2.2. *Differential operators.* In this section, we introduce the differential operators which represent the analogs of the Malliavin derivative and of the Skorohod integral.

*Simple functionals.* A random variable $F$ is called a *simple functional* if there exists some $n$ and some $\mathcal{G} \times \mathcal{B}(\mathbb{R}^n)$-measurable function $f : \Omega \times \mathbb{R}^n \to \mathbb{R}$ such that

$$F = f(\omega, V_1, \ldots, V_n).$$

We denote by $\mathcal{S}_{(n,k)}$ the space of simple functionals such that $f \in \mathcal{C}_{n,k}$, and $\mathcal{S}_{(n,k)}(A)$ will denote the space of simple functionals such that $f \in \mathcal{C}_{n,k}(A)$.

We will use the notation $\partial_{V_i} F := \frac{\partial f}{\partial x_i}(\omega, V_1, \ldots, V_n)$, $i = 1, \ldots, n$.

*Simple processes.* A simple process of length $n$ is a sequence of random variables $U = (U_i)_{i \le n}$ such that

$$U_i(\omega) = u_i(\omega, V_1(\omega), \ldots, V_n(\omega)),$$

where $u_i : \Omega \times \mathbb{R}^n \to \mathbb{R}$, $i \in \mathbb{N}$, are $\mathcal{G} \times \mathcal{B}(\mathbb{R}^n)$-measurable functions. We denote by $\mathcal{P}_{(n,k)}$ [resp., $\mathcal{P}_{(n,k)}(A)$] the space of simple processes of length $n$ such that $u_i \in \mathcal{C}_{n,k}$, $i = 1, \ldots, n$ [resp., $u_i \in \mathcal{C}_{n,k}(A)$, $i = 1, \ldots, n$]. Note that if $U \in \mathcal{P}_{(n,k)}$, then $U_i \in \mathcal{S}_{(n,k)}$ and if $U \in \mathcal{P}_{(n,k)}(A)$, then $U_i \in \mathcal{S}_{(n,k)}(A)$.

On the space of simple processes, we consider the scalar product

$$\langle U, V \rangle_\pi := \sum_{i=1}^n \pi_i(\omega, V_i) U_i(\omega) V_i(\omega).$$

We now define the differential operators which appear in Malliavin's calculus.

*The Malliavin derivative.* $D : \mathcal{S}_{(n,1)} \to \mathcal{P}_{(n,0)}$: if $F = f(\omega, V_1, \ldots, V_n)$, then

$$D_i F := \frac{\partial f}{\partial x_i}(\omega, V_1(\omega), \ldots, V_n(\omega)) \mathbf{1}_{B_i(\omega)}(V_i),$$
$$DF = (D_i F)_{i \le n} \in \mathcal{P}_{(n,0)}.$$

*The Malliavin covariance matrix.* Given $F=(F^1,\dots,F^d)$, $F^i=f^i(\omega,V_1,\dots,V_n)\in\mathcal{S}_{(n,1)}$, the Malliavin covariance matrix is

$$\sigma_F^{ij}=\langle DF^i,DF^j\rangle_\pi=\sum_{p=1}^n \pi_p(\omega,V_p)\,\partial_p f^i\,\partial_p f^j(\omega,V_1,\dots,V_n).$$

This is a symmetric, positive-definite matrix.

*The Skorohod integral.* We define $\delta:\mathcal{P}_{(n,1)}\to\mathcal{S}_{(n,0)}$: for $U=(U_i)_{1\le i\le n}$ such that $U_i(\omega)=u_i(\omega,V_1,\dots,V_n)$, we define

$$\delta_i(U):=-\biggl(\frac{\partial}{\partial x_i}(\pi_i u_i)+(\pi_i u_i)\,\partial\ln p_i\biggr)(\omega,V_1,\dots,V_n),$$

$$\delta(U):=\sum_{i=1}^n\delta_i(U).$$

*The border term operator.* For $F=f(\omega,V_1,\dots,V_n)\in\mathcal{S}_{(n,0)}$ and $U=(u_i(\omega,V_1,\dots,V_n))_{i=1,\dots,n}\in\mathcal{P}_{(n,0)}$, we define

$$\begin{aligned}
[F,U]_\pi &= \sum_{i=1}^n \Gamma_i(I_i(f\times u_i)\times\pi_i\times p_i)\\
&=\sum_{i=1}^n\sum_{j=1}^{k_i}((f\times u_i)(\omega,V_1,\dots,V_{j-1},t_i^j-,V_{j+1},\dots,V_n)(\pi_i p_i)(\omega,t_i^j-)\\
&\qquad -(f\times u_i)(\omega,V_1,\dots,V_{j-1},\\
&\qquad\qquad t_i^j+,V_{j+1},\dots,V_n)(\pi_i p_i)(\omega,t_i^j+))\\
&\quad+\sum_{i=1}^n(f\times u_i)(\omega,V_1,\dots,V_{j-1},b_i-,V_{j+1},\dots,V_n)(\pi_i p_i)(\omega,b_i-)\\
&\quad-\sum_{i=1}^n(f\times u_i)(\omega,V_1,\dots,V_{j-1},a_i+,V_{j+1},\dots,V_n)(\pi_i p_i)(\omega,a_i+).
\end{aligned}$$

REMARK 2.1. If we choose $\pi_i$ such that

$$\begin{aligned}
&\pi_i(\omega,t_i^j+)=\pi_i(\omega,t_i^j-)=0, \qquad i=1,\dots,n,j=1,\dots,k_i,\\
&\pi_i(\omega,a_i+)=\pi_i(\omega,b_i-)=0, \qquad i=1,\dots,n,
\end{aligned}\tag{2.3}$$

then $[F,U]_\pi=0$ for every $F\in\mathcal{S}_{(n,1)}$ and $U\in\mathcal{P}_{(n,1)}$. So there will be no border terms in the duality formula or in the integration by parts formula. This is one possible reason for introducing the weights. The other one concerns renormalization.

In our framework the duality between $\delta$ and $D$ is given by the following proposition.

PROPOSITION 2.1. *Let $F \in \mathcal{S}_{(n,1)}$ and $U \in \mathcal{P}_{(n,1)}$. Suppose that for every $i = 1, \dots, n$,*

$$\mathrm{E}(|F\delta_i(U)|\mathbf{1}_A) + \mathrm{E}(\pi_i(\omega, V_i)|D_iF \times U_i|\mathbf{1}_A) < \infty. \tag{2.4}$$

*Then $\mathrm{E}(|[F,U]_\pi|\mathbf{1}_A) < \infty$ and*

$$\mathrm{E}(\langle DF, U\rangle_\pi \mathbf{1}_A) = \mathrm{E}(F\delta(U)\mathbf{1}_A) + \mathrm{E}([F,U]_\pi \mathbf{1}_A). \tag{2.5}$$

*If (2.3) holds true, then*

$$\mathrm{E}(\langle DF, U\rangle_\pi \mathbf{1}_A) = \mathrm{E}(F\delta(U)\mathbf{1}_A).$$

PROOF. Since $\pi_i = 0$ on $(a_i, b_i)^c$, we have

$$\begin{aligned}
&\mathrm{E}(\langle DF, U\rangle_\pi \mathbf{1}_A) \\
&\quad = \mathrm{E}\left(\sum_{i=1}^n \mathrm{E}(\pi_i(\omega, V_i) D_i F \times U_i \mid \mathcal{G}_i)\mathbf{1}_A\right) \\
&\quad = \mathrm{E}\left(\mathbf{1}_A \sum_{i=1}^n \int_{a_i}^{b_i} (\pi_i u_i\, \partial_i f)(\omega, V_1, \dots, V_{i-1}, y, V_{i+1}, \dots, V_n) p_i(\omega, y)\, dy\right).
\end{aligned}$$

Using integration by parts [see (2.2)], we obtain

$$\begin{aligned}
&\int_{a_i}^{b_i} \partial_i f \times (\pi_i u_i) \times p_i \\
&\quad = \sum_{j=0}^{k_i} \int_{(t_i^j, t_{i+1}^j)} \partial_i f \times (\pi_i u_i) \times p_i \\
&\quad = \Gamma_i(I_i(f \times u_i)\pi_i p_i) - \sum_{j=0}^{k_i} \int_{(t_i^j, t_{i+1}^j)} f \times (\partial_i(\pi_i u_i) \times p_i + (\pi_i u_i) \times \partial p_i) \\
&\quad = \Gamma_i(I_i(f \times u_i)\pi_i p_i) - \int_{a_i}^{b_i} f \times (\partial_i(\pi_i u_i) + \pi_i u_i\, \partial \ln p_i) \times p_i.
\end{aligned}$$

By (2.4), we have

$$\int_{\mathbb{R}} (|u_i\, \partial_i f|\pi_i p_i)(\omega, V_1, \dots, V_{i-1}, y, V_{i+1}, \dots, V_n)\, dy < \infty,$$

$$\int_{\mathbb{R}} (|f(\partial_i(\pi_i u_i) + \pi_i u_i\, \partial \ln p_i)| \times p_i)(\omega, V_1, \dots, V_{i-1}, y, V_{i+1}, \dots, V_n)\, dy < \infty,$$

for almost all $\omega \in A$. Thus, the above integrals make sense. Since $\Gamma_i(I_i(f \times u_i)\pi_i p_i)$ is the sum of these two integrals, we also obtain $\mathrm{E}(|\Gamma_i(I_i(f \times u_i)\pi_i p_i)| \mathbf{1}_A) < \infty$ so that $\mathrm{E}(|[F,U]_\pi| \mathbf{1}_A) < \infty$.

Using the definition of $p_i$, we return to expectations and obtain

$$\int_{a_i}^{b_i} (\pi_i u_i \, \partial_i f)(\omega, V_1, \ldots, V_{i-1}, y, V_{i+1}, \ldots, V_n) p_i(\omega, y)\, dy$$
$$= \mathrm{E}(F\delta_i(U) \mid \mathcal{G}_i) + \Gamma_i(I_i(f \times u_i)\pi_i p_i).$$

One sums over $i$ and the proof is complete. $\square$

Corollary 2.1. *Let $Q \in \mathcal{S}_{(n,1)}(A)$ satisfy*

$$\mathrm{E}(\mathbf{1}_A(|\pi_i Q| + |\partial_{V_i}(\pi_i Q)|)^{1+\eta}) < \infty, \qquad i = 1, \ldots, n, \tag{2.6}$$

*for some $\eta > 0$. Then for every $F \in \mathcal{S}_{(n,1)}(A)$, $U \in \mathcal{P}_{(n,1)}(A)$, one has*

$$\mathrm{E}(Q\langle DF, U\rangle_\pi \mathbf{1}_A) = \mathrm{E}(F\delta(QU)\mathbf{1}_A) + \mathrm{E}([F, QU]_\pi \mathbf{1}_A). \tag{2.7}$$

Proof. We need only check that $F$ and $\widetilde{U} = QU$ satisfy (2.4). We have

$$|\delta_i(QU)| \le |\partial_{V_i}(\pi_i Q)||U_i| + |\pi_i Q|(|\partial_{V_i} U_i| + |U_i||\partial \ln p_i|).$$

Since $U \in \mathcal{P}_{(n,1)}(A)$, one has $U_i, \partial_{V_i} U_i \in L_{(\infty)}(A)$ and by Hypothesis 2.2, $\partial \ln p_i \in L_{(\infty)}(A)$. So, using (2.6), we have $\delta_i(QU) \in L_{(1+)}(A)$ and since $F \in L_{(\infty)}(A)$, we obtain $\mathrm{E}(F\delta_i(QU)|) < \infty$.

We have $D_i F$, $U_i \in L_{(\infty)}(A)$ and $\pi_i Q \in L_{(1+)}(A)$, so $\mathrm{E}(\pi_i |D_i F \times (QU_i)|) < \infty$. $\square$

*The Ornstein–Uhlenbeck operator.* We now introduce $L := \delta(D) : \mathcal{S}_{(n,2)} \to \mathcal{S}_{(n,0)}$:

$$LF := -\sum_{i=1}^n (\partial_i(\pi_i\, \partial_i f) + \pi_i\, \partial_i f \partial \ln p_i)(\omega, V_1, \ldots, V_n)$$
$$= -\sum_{i=1}^n ((\pi_i' + \pi_i\, \partial \ln p_i)\partial_i f + \pi_i\, \partial_i^2 f)(\omega, V_1, \ldots, V_n).$$

As an immediate consequence of the duality relation (2.5), we obtain the following.

Lemma 2.1. *Let $F, G \in \mathcal{S}_{(n,2)}$ and $A \in \mathcal{G}$. Suppose that for every $i = 1, \ldots, n$,*

$$\mathrm{E}[(|FL_i G| + |GL_i F| + \pi_i |D_i F \times D_i G|)\mathbf{1}_A] < \infty.$$

*Then $\mathrm{E}(|[F, DG]_\pi| \mathbf{1}_A) < \infty$, $\mathrm{E}(|[G, DF]_\pi| \mathbf{1}_A) < \infty$ and*

$$\mathrm{E}(FLG\mathbf{1}_A) + \mathrm{E}([F, DG]_\pi \mathbf{1}_A) = \mathrm{E}(\langle DF, DG\rangle_\pi \mathbf{1}_A)$$
$$= \mathrm{E}(GLF\mathbf{1}_A) + \mathrm{E}([G, DF]_\pi \mathbf{1}_A).$$

We denote by $\mathcal{C}_p^k(\mathbb{R}^d)$ the space of the functions $\phi:\mathbb{R}^d\to\mathbb{R}$ which are $k$-times differentiable and such that $\phi$ and its derivatives of order less then or equal to $k$ have polynomial growth. Standard differential calculus gives the following chain rules.

LEMMA 2.2. (i) *Let* $\phi\in\mathcal{C}_p^1(\mathbb{R}^d)$ *and* $F=(F^1,\dots,F^d)$, $F^i\in\mathcal{S}_{(n,1)}(A)$. *Then* $\phi(F)\in\mathcal{S}_{(n,1)}(A)$ *and*

$$D\phi(F)=\sum_{k=1}^d \partial_k\phi(F)DF^k. \tag{2.8}$$

(ii) *If* $\phi\in\mathcal{C}_p^2(\mathbb{R}^d)$ *and* $F^i\in\mathcal{S}_{(n,2)}(A)$, *then* $\phi(F)\in\mathcal{S}_{(n,2)}(A)$ *and*

$$L\phi(F)=\sum_{k=1}^d \partial_k\phi(F)LF^k-\sum_{k,p=1}^d \partial_{k,p}^2\phi(F)\langle DF^k,DF^p\rangle_\pi.$$

(iii) *Let* $F\in\mathcal{S}_{(n,1)}(A)$ *and* $U\in\mathcal{P}_{(n,1)}(A)$. *Then* $FU\in\mathcal{P}_{(n,1)}(A)$ *and*

$$\delta(FU)=F\delta(U)-\langle DF,U\rangle_\pi.$$

*In particular, if* $F\in\mathcal{S}_{(n,1)}(A)$ *and* $G\in\mathcal{S}_{(n,2)}(A)$, *then* $FDG\in\mathcal{P}_{(n,1)}(A)$ *and*

$$\delta(FDG)=FLG-\langle DF,DG\rangle_\pi. \tag{2.9}$$

REMARK 2.2. Let us define $L^2_{\pi,n}(A)$ to be the closure of $\mathcal{P}_{(n,0)}$ with respect to the norm associated with the scalar product $\langle U,V\rangle=\mathrm{E}(\langle U,V\rangle_\pi)$. If $[F,U]_\pi$ is not null, then the operator $D:\mathcal{S}_{(n,1)}\subset L^2(\Omega)\to\mathcal{P}_{(n,0)}\subset L^2_{\pi,n}(A)$ is not closable.

Suppose, for example, that $V_1$ is exponentially distributed and $V_i$, $i=2,\dots,n$, are arbitrary and independent of $V_1$. We take $\pi_1=1$ and $\pi_i=0$, $i=2,\dots,n$. So we define our calculus with respect to $V_1$ only. In this case, $a_1=0,b_1=\infty$ and there are no points $t_i^j$. Now take $F_n=f_n(V_1)$ with $f_n(x)=1-nx$ for $0<x<1/n$ and $f_n(x)=0$ for $x\geq 1/n$. Also, take $u_1(x)=1-x$ for $0<x<1$ and $u_1(x)=0$ for $x\geq 1$ and write the duality formula $\mathrm{E}(\langle DF_n,U\rangle_\pi)=\mathrm{E}(F_n\delta(U))+\mathrm{E}([F_n,U]_\pi)$. Since $[F_n,U]_\pi=1$ and $F_n\to 0$ in $L^2(\Omega)$, we obtain $\lim_{n\to\infty}\mathrm{E}(\langle DF_n,U\rangle_\pi)=1$ and so $DF_n\nrightarrow 0$ in $L^2_{\pi,n}(A)$. This proves that $D$ is not closable.

But if $[F,U]_\pi=0$ for every $F$, $U$ [this happens, e.g., if we choose $\pi_i$ to satisfy (2.3)], then the duality formula (2.5) guarantees that $D$ and $\delta$ are closable. But we will remain at the level of simple functionals and will not discuss the extension to the infinite-dimensional setting.

REMARK 2.3. The above differential operators and the duality formula (2.5) represent abstract versions of the operators introduced in Malliavin calculus and of the duality formula used there. In order to see this, we consider the simple example of the Euler scheme for a diffusion process, corresponding to the time grid $0 = s_0 < s_1 < \cdots < s_n = s$. This is a simple functional depending on the increments of the Brownian motion $B$, that is, $V_i = B(s_i) - B(s_{i-1})$, $i = 1, \ldots, n$. The variables on which the calculus is based are independent Gaussian variables. It follows that $p_i(\omega, y) = (2\pi(s_i - s_{i-1}))^{-1/2} \times \exp(-y^2/2(s_i - s_{i-1}))$. Since $p_i$ is smooth on the whole of $\mathbb{R}$ and has null limit at infinity, there will be no border terms, so we take $a_i = -\infty, b_i = \infty$ and $k_i = 0$. If $F = f(V_1, \ldots, V_n)$, then $D_i F = \partial_i f(V_1, \ldots, V_n) = \overline{D}_s F \mathbf{1}_{[s_{i-1}, s_i)}(s)$, where $\overline{D}_s F$ is the standard Malliavin derivative. We take $\pi_i = s_i - s_{i-1}$ so that

$$\langle DF, DG \rangle_\pi = \sum_{i=1}^n \pi_i D_i F D_i G = \int_0^s \overline{D}_u F \overline{D}_u G \, du.$$

We note that here the weights are used in order to obtain the Lebesgue measure. Moreover, we have $\partial_y \ln p_i(y) = -y/(s_i - s_{i-1})$, so

$$\delta_i(U) = -\sum_{i=1}^n (\partial_i u_i(V_1, \ldots, V_n)(s_i - s_{i-1}) - u_i(V_1, \ldots, V_n) V_i).$$

So we find out the standard Malliavin calculus.

REMARK 2.4. If $[F, G]_\pi = 0$, the calculus presented here fits into the framework introduced by Bouleau in [5]: in the notation there, the bilinear form $(F, G) \to \langle DF, DG \rangle_\pi$ leads to an error structure. A variety of examples and applications of these structures is discussed. That framework mainly focuses on the error calculus, but examples of applications to sensitivity computations are also given and an integration by parts formula is derived. This works well in the particular case of a one-dimensional functional. Moreover, the differential calculus is based on a single noise $V_i$ as in Corollary 2.2 below (so the weights $\pi_i$ do not come into the nondegeneracy condition). In a more general framework, the nondegeracy condition involves the weights $\pi_i$, $i \in \mathbb{N}$, and a more detailed analysis must be undertaken (see the following section).

2.3. *The integration by parts formula.* We consider $F = (F^1, \ldots, F^d) \in \mathcal{S}^d_{(n,1)}(A)$ and define

$$\Theta_F(A) := \{G = \sigma_F \times Q : Q \in \mathcal{S}^d_{(n,1)}(A), Q_i \text{ satisfy } (2.6)\}.$$

We think of $G \in \Theta_F(A)$ as a random direction in which $F$ is nondegenerate (in Malliavin's sense).

The basic integration by parts formula is given in the following theorem.

THEOREM 2.1. *Let* $F=(F^1,\dots,F^d)\in\mathcal{S}^d_{(n,2)}(A)$ *and* $G\in\Theta_F(A)$, $G=\sigma_F\times Q$. *Then*

$$\delta\Bigg(\sum_{i=1}^d Q^i DF^i\Bigg),\Bigg[\phi(F),\sum_{i=1}^d Q^i DF^i\Bigg]_\pi\in L_{(1+)}(A)$$

*and for every* $\phi\in\mathcal{C}^1_p(\mathbb{R}^d)$, *one has*

$$\begin{aligned}\mathrm{E}(\langle\nabla\phi(F),G\rangle\mathbf{1}_A)&=\mathrm{E}\Bigg(\phi(F)\delta\Bigg(\sum_{i=1}^d Q^i DF^i\Bigg)\mathbf{1}_A\Bigg)\\&\quad+\mathrm{E}\Bigg(\Bigg[\phi(F),\sum_{i=1}^d Q^i DF^i\Bigg]_\pi\mathbf{1}_A\Bigg).\end{aligned}\tag{2.10}$$

PROOF. Using (2.8),

$$\langle D\phi(F),DF^i\rangle_\pi=\sum_{j=1}^d\partial_j\phi(F)\langle DF^j,DF^i\rangle_\pi=\sum_{j=1}^d\partial_j\phi(F)\sigma_F^{ij}.$$

Since $G=\sigma_F\times Q$, we obtain

$$\begin{aligned}\langle\nabla\phi(F),G\rangle&=\sum_{j=1}^d\partial_j\phi(F)G^j=\sum_{j=1}^d\partial_j\phi(F)\sum_{i=1}^d Q^i\sigma_F^{ij}\\&=\sum_{i=1}^d Q^i\sum_{j=1}^d\partial_j\phi(F)\sigma_F^{ij}=\sum_{i=1}^d Q^i\langle D\phi(F),DF^i\rangle_\pi.\end{aligned}$$

Note that $\phi(F)\in\mathcal{S}_{(n,1)}(A)$ and $DF^i\in\mathcal{P}_{(n,1)}(A)$. Since the $Q_i$ satisfy (2.6), one may use the duality formula (2.7) and thereby obtain (2.10). □

We now give a nondegeneracy condition on $\sigma_F$ which guarantees that all of the directions are nondegenerate for $F$.

We assume that $\det\sigma_F\neq 0$ on $A$ and define $\gamma_F=\sigma_F^{-1}$. We also assume that $\pi_l(\det\gamma_F)^2,\pi_l'\det\gamma_F,\pi_l\pi_l'(\det\gamma_F)^2\in L_{(1+)}(A)$ for every $l=1,\dots,n$. This means that there exists $\eta>0$ such that

$$\mathrm{E}[\mathbf{1}_A(|\pi_l|(\det\gamma_F)^2+|\pi_l'|(\det\gamma_F+|\pi_l|(\det\gamma_F)^2))^{1+\eta}]<\infty.\tag{2.11}$$

LEMMA 2.3. *Assume that* (2.11) *holds true and that* $F\in\mathcal{S}^d_{(n,2)}(A)$. *Then* $\mathcal{S}^d_{(n,1)}(A)\subseteq\Theta_F(A)$.

PROOF. Let $G \in S^d_{(n,1)}(A)$. Then $G = \sigma_F \times Q$ with $Q = \gamma_F \times G$. We write $\gamma_F^{ij} = \widehat{\sigma}_F^{ij} \times \det\gamma_F$, where $\widehat{\sigma}_F^{ij}$ is the algebraic complement. It follows that $Q^i = \det\gamma_F \times S^i$ with $S^i = \sum_{j=1}^d G^j \widehat{\sigma}_F^{ij}$.

Let us check that (2.6) holds true for $Q^i$, $i = 1, \ldots, d$. Since $\pi_l \in L_{(\infty)}(A)$ and $D_l F^i \in L_{(\infty)}(A)$, one has $\widehat{\sigma}_F^{ij}$, $\det\sigma_F \in L_{(\infty)}(A)$ and since $G^j \in L_{(\infty)}(A)$, we have $S^i \in L_{(\infty)}(A)$. Moreover, by (2.11), $\pi_l \det\gamma_F \in L_{(1+)}(A)$, so $\pi_l Q^i = (\pi_l \det\gamma_F) S^i \in L_{(1+)}(A)$.

We now check that $D_l(\pi_l Q^i) \in L_{(1+)}(A)$. We write

$$D_l \sigma_F^{ij} = \pi_l' D_l F^i D_l F^j + \sum_{k=1}^n \pi_k D_l(D_k F^i D_k F^j).$$

Since $F \in S^d_{(n,2)}(A)$, we have $D_l F^i D_l F^j$, $D_l(D_k F^i D_k F^j) \in L_{(\infty)}(A)$ and consequently $D_l \sigma_F^{ij} = \theta_1 + \theta_2 \pi_l'$ with $\theta_1$, $\theta_2 \in L_{(\infty)}(A)$. Then $D_l(\det\sigma_F) = \mu + \nu\pi_l'$ and $D_l S^i = \mu_i + \nu_i \pi_l'$ with $\mu$, $\nu$, $\mu_i$, $\nu_i \in L_{(\infty)}(A)$.

Using (2.11), we obtain

$$\begin{aligned} D_l(\pi_l Q^i) &= \pi_l' \det\gamma_F S^i - \pi_l (\det\gamma_F)^2 D_l(\det\sigma_F) S^i + \pi_l \det\gamma_F D_l S^i \\ &= \pi_l' \det\gamma_F S^i - \pi_l (\det\gamma_F)^2 (\mu + \nu\pi_l') S^i + \pi_l \det\gamma_F (\mu_i + \nu_i \pi_l') \\ &\in L_{(1+)}(A) \end{aligned}$$

and the proof is complete. □

As a consequence, we obtain the following.

THEOREM 2.2. *Let $F = (F^1, \ldots, F^d) \in \mathcal{S}^d_{(n,2)}(A)$ and $G \in \mathcal{S}_{(n,1)}(A)$. Suppose that (2.11) holds true. Then*

$$\delta\left(G \sum_{j=1}^d \gamma_F^{ji} DF^j\right), \left[\phi(F), G \sum_{j=1}^d \gamma_F^{ji} DF^j\right]_\pi \in L_{(1+)}(A)$$

*and for every $\phi \in \mathcal{C}^1_p(\mathbb{R}^d)$, one has*

$$\mathrm{E}(\partial_i \phi(F) G \mathbf{1}_A) = \mathrm{E}\left[\phi(F)\delta\left(G \sum_{j=1}^d \gamma_F^{ji} DF^j\right)\mathbf{1}_A\right] + \mathrm{E}\left(\left[\phi(F), G \sum_{j=1}^d \gamma_F^{ji} DF^j\right]_\pi \mathbf{1}_A\right)$$

*for every $i = 1, \ldots, d$.*

*Suppose that $\pi_l$, $l = 1, \ldots, n$, satisfy (2.3). Then*

$$\mathrm{E}(\partial_i \phi(F) G \mathbf{1}_A) = \mathrm{E}(\phi(F) H_i(F, G) \mathbf{1}_A)$$

*with*

$$H_i(F,G) = \delta\left(G\sum_{j=1}^{d}\gamma_F^{ji}DF^j\right) = \sum_{j=1}^{d}(G\gamma_F^{ji}LF^j - \langle D(G\gamma_F^{ji}), DF^j\rangle_\pi).$$

PROOF. We take $\widetilde{G} = (0,\dots,0,G,0,\dots,0)$ with $G$ occupying the $i$th place so that $\partial_i\phi(F)G = \langle\nabla\phi(F),\widetilde{G}\rangle$. In view of Lemma 2.3, $\widetilde{G} \in \Theta_F(A)$ and $\widetilde{G} = \sigma_F \times Q$, with $Q^j = G\gamma_F^{ji}$. One then employs Theorem 2.1. In order to obtain the second equality in the expression for $H_i(F,G)$, one employs (2.9). □

There is one particular situation in which the nondegeneracy condition (2.11) does not involve the weights—when if $F$ is one-dimensional and the integration by parts formula is based on a single random variable $V_i$. We then have the following corollary.

COROLLARY 2.2. *Let $F = f(V_1,\dots,V_n) \in \mathcal{S}_{(n,2)}(A)$ and $G \in \mathcal{S}_{(n,1)}(A)$. Suppose that there exists some $l \in \{1,\dots,n\}$ such that*

$$\mathrm{E}[\mathbf{1}_A(D_lF)^{-6(1+\eta)}] < \infty \tag{2.12}$$

*for some $\eta > 0$. Consider the weights $\pi_i = 0$ for $i \neq l$ and let $\pi_l$ be an arbitrary function which satisfies $\pi_l \in L_{(\infty)}(A)$ and $\pi_l' \in L_{(1+)}(A)$. Then $\delta(G\gamma_F DF)$, $[\phi(F), G\gamma_F DF]_\pi \in L_{(1+)}(A)$ and for every $\phi \in \mathcal{C}_p^1(\mathbb{R})$, one has*

$$\mathrm{E}(\phi'(F)G\mathbf{1}_A) = \mathrm{E}(\phi(F)\delta(G\gamma_F DF)\mathbf{1}_A) + \mathrm{E}([\phi(F), G\gamma_F DF]_\pi\mathbf{1}_A). \tag{2.13}$$

PROOF. Note that $\sigma_F = \pi_l(V_l)|D_lF|^2$. We return to the proof of Theorem 2.1 and write $G = Q\sigma_F$ with

$$Q = \begin{cases} \dfrac{G}{\pi_l(V_l)|D_lF|^2}, & \text{if } \pi_l(V_l)|D_lF|^2 \neq 0, \\ 0, & \text{if } \pi_l(V_l)|D_lF|^2 = 0. \end{cases}$$

Then $\pi_l(V_l)Q = g(V_1,\dots,V_n)/|D_lF|^2$ and, as a consequence of the hypothesis (2.12), one has $\pi_l(V_l)Q$, $\partial_{V_i}(\pi(V_l)Q) \in L_{(1+)}(A)$, $i = 1,\dots,n$. So we may use the duality relation to conclude the proof. □

## 3. Pure jump diffusions.

In this section, we will use the integration by parts formula presented in the previous section for a pure jump diffusion $(S_t)_{t\geq 0}$. We will use the notation from [14]. We consider a Poisson point measure $N(dt,da)$ on $\mathbb{R}$ with positive and finite intensity measure $\mu(da)\times dt$, that is, $\mathrm{E}(N([0,t]\times A)) = \mu(A)t$. We denote by $J_t$ the counting process, that is, $J_t := N([0,t]\times\mathbb{R})$ and we denote by $T_i$, $i \in \mathbb{N}$, the jump times of $J_t$.

We represent the above Poisson point measure by means of a sequence $\Delta_i$, $i \in \mathbb{N}$, of independent random variables with law $\nu(da) = \mu(\mathbb{R})^{-1} \times \mu(da)$. This means that $N([0,t] \times A) = \mathit{card}\{T_i \leq t : \Delta_i \in A\}$.

We look at the solution $S_t$ of the equation

$$
\begin{aligned}
S_t &= x + \sum_{i=1}^{J_t} c(T_i, \Delta_i, S_{T_i^-}) + \int_0^t g(r, S_r)\, dr, \\
&= x + \int_0^t \int_{\mathbb{R}} c(s, a, S_{s^-})\, dN(s,a) + \int_0^t g(r, S_r)\, dr, \qquad 0 \leq t \leq T.
\end{aligned}
\tag{3.1}
$$

We will work under the following hypothesis.

HYPOTHESIS 3.1. The functions $(t,x) \to c(t,a,x)$ and $x \to g(t,x)$ are twice differentiable and have bounded derivatives of first and second order. Moreover, we assume that they have linear growth with respect to $x$, uniformly with respect to $t$ and $a$, that is, $|c(t,a,x)| + |g(t,x)| \leq K(1+|x|)$.

On each set $\{J_t = n\}$, $S_t$ is a simple functional of $\Delta_1, \ldots, \Delta_n$ and $T_1, \ldots, T_n$. In the first subsection, we present the deterministic calculus which permits us to compute the Malliavin derivatives and in the following two subsections, we give the integration by parts formula with respect to the amplitude of the jumps and with respect to the jump times, separately. Finally, in the fourth subsection, we briefly present the mixed calculus with respect to both.

We shall remain in the one-dimensional case (i.e., $S_t \in \mathbb{R}$) because the multi-dimensional case is more involved from a technical point of view. Our purpose is to illustrate the way in which the integration by parts formula works for Poisson point measures and to emphasize the specific difficulties. The heavy techniques related to the multi-dimensional case would obscure these specific points, but the machinery works just as well in this case.

3.1. *The deterministic equation.* We fix some deterministic $0 = u_0 < u_1 < \cdots < u_n < T$ and define $u = (u_1, \ldots, u_n)$. We also fix $a = (a_1, \ldots, a_n) \in \mathbb{R}^n$. To these fixed numbers, we associate the deterministic equation

$$
s_t = x + \sum_{i=1}^{J_t(u)} c(u_i, a_i, s_{u_i^-}) + \int_0^t g(r, s_r)\, dr, \qquad 0 \leq t \leq T, \tag{3.2}
$$

where $J_t(u) = k$ if $u_k \leq t < u_{k+1}$. We denote by $s_t(u,a)$, or simply by $s_t$, the solution of this equation. This is the deterministic counterpart of our stochastic equation. On the set $\{J_t = n\}$, the solution $S_t$ of (3.1) is represented as $S_t = s_t(T_1, \ldots, T_n, \Delta_1, \ldots, \Delta_n)$.

In order to solve this equation, we introduce the flow $\Phi = \Phi_u(t,x)$, $0 \le u \le t$, $x \in \mathbb{R}$, which solves the ordinary integral equation

$$\Phi_u(t,x) = x + \int_u^t g(r, \Phi_u(r,x))\, dr, \qquad t \ge u. \tag{3.3}$$

The solution $s$ of the equation (3.2) is given by

$$\begin{aligned} s_0 &= x, \\ s_t &= \Phi_{u_i}(t, s_{u_i}) \qquad \text{for } u_i \le t < u_{i+1}, \\ s_{u_{i+1}} &= s_{u_{i+1}^-} + c(u_{i+1}, a_{i+1}, s_{u_{i+1}^-}) \\ &= \Phi_{u_i}(u_{i+1}, s_{u_i}) + c(u_{i+1}, a_{i+1}, \Phi_{u_i}(u_{i+1}, s_{u_i})). \end{aligned} \tag{3.4}$$

Our aim is to compute the derivatives of $s$ with respect to $u_j$, $a_j$, $j = 1, \dots, n$.

We first introduce some notation. We define

$$e_{u,t}(x) := \exp\biggl(\int_u^t \partial_x g(r, \Phi_u(r,x))\, dr\biggr).$$

Since $\Phi_{u_i}(r, s_{u_i}) = s_r$ for $u_i \le r < u_{i+1}$, we have

$$e_{u_i,t}(s_{u_i}) = \exp\biggl(\int_{u_i}^t \partial_x g(r, s_r)\, dr\biggr) \qquad \text{for } u_i \le t < u_{i+1}.$$

Since

$$\partial_x \Phi_u(t,x) = 1 + \int_u^t \partial_x g(r, \Phi_u(r,x))\, \partial_x \Phi_u(r,x)\, dr,$$

it follows that

$$\partial_x \Phi_u(t,x) = e_{u,t}(x)$$

and since

$$\partial_u \Phi_u(t,x) = -g(u,x) + \int_u^t \partial_x g(r, \Phi_u(r,x))\, \partial_u \Phi_u(r,x)\, dr,$$

we have

$$\partial_u \Phi_u(t,x) = -g(u,x) e_{u,t}(x).$$

Finally, we define

$$q(t,\alpha,x) := (\partial_t c + g\, \partial_x c)(t,\alpha,x) + g(t,x) - g(t, x + c(t,\alpha,x)).$$

LEMMA 3.1. *Suppose that Hypothesis 3.1 holds true. Then $s_t(u,a)$ is twice differentiable with respect to $u_j, j=1,\ldots,n$, and with respect to $a_j, j=1,\ldots,n$, and we have the following explicit expressions for the derivatives.*

A. Derivatives with respect to $u_j$. *For $t<u_j$, $\partial_{u_j}s_t(u,a)=0$. Moreover,*

$$\partial_{u_j}s_{u_j-}=g(u_j,s_{u_j-}),$$
$$\partial_{u_j}s_{u_j}=(\partial_t c+g(1+\partial_x c))(u_j,a_j,s_{u_j-}).$$

*For $u_j<t<u_{j+1}$,*

$$\begin{aligned}\partial_{u_j}s_t&=q(u_j,a_j,s_{u_j-})e_{u_j,t}(s_{u_j}),\\ \partial_{u_j}s_{u_{j+1}-}&=q(u_j,a_j,s_{u_j-})e_{u_j,u_{j+1}}(s_{u_j}),\\ \partial_{u_j}s_{u_{j+1}}&=q(u_j,a_j,s_{u_j-})(1+\partial_x c(u_{j+1},a_{j+1},s_{u_{j+1}-}))e_{u_j,u_{j+1}}(s_{u_j}).\end{aligned}\tag{3.5}$$

*Finally, for $p\geq j+1$ and $u_p\leq t<u_{p+1}$, we have the recurrence relations*

$$\begin{aligned}\partial_{u_j}s_t&=e_{u_p,t}(s_{u_p})\,\partial_{u_j}s_{u_p},\\ \partial_{u_j}s_{u_{p+1}}&=(1+\partial_x c(u_{p+1},a_{p+1},s_{u_{p+1}-}))e_{u_p,u_{p+1}}(s_{u_p})\,\partial_{u_j}s_{u_p}.\end{aligned}\tag{3.6}$$

*Define $T(f):=\partial_t f+g\,\partial_x f$. The second order derivatives are given by*

$$\partial^2_{u_j}s_{u_j-}=T(g)(u_j,a_j,s_{u_j-}),$$
$$\partial^2_{u_j}s_{u_j}=T(\partial_t c+g(1+\partial_x c))(u_j,a_j,s_{u_j-}).$$

*Define*

$$\begin{aligned}\rho_j(t)&=\partial_{u_j}e_{u_j,t}(s_{u_j})\\ &=e_{u_j,t}(s_{u_j})\biggl(-\partial_x g(u_j,s_{u_j})+q(u_j,a_j,s_{u_j-})\int_{u_j}^t\partial_x^2 g(r,s_r)e_{u_j,r}(s_{u_j})\,dr\biggr).\end{aligned}$$

*Then for $u_j<t<u_{j+1}$,*

$$\partial^2_{u_j}s_t(u,a)=T(q)(u_j,a_j,s_{u_j-}(u,a))e_{u_j,t}(s_{u_j})+q(u_j,a_j,s_{u_j-}(u,a))\rho_j(t)$$

*and*

$$\begin{aligned}\partial^2_{u_j}s_{u_{j+1}}&=T(q)(u_j,a_j,s_{u_j-})(1+\partial_x c)(u_{j+1},a_{j+1},s_{u_{j+1}-})e_{u_j,u_{j+1}}(s_{u_j})\\ &\quad+q^2(u_j,a_j,s_{u_j-})\,\partial_x^2 c(u_{j+1},a_{j+1},s_{u_{j+1}-})e^2_{u_j,u_{j+1}}(s_{u_j})\\ &\quad+q(u_j,a_j,s_{u_j-})(1+\partial_x c)(u_{j+1},a_{j+1},s_{u_{j+1}-})\rho_j(u_j).\end{aligned}$$

*For $p\geq j+1$, we define*

$$\rho_{j,p}(t)=\partial_{u_j}e_{u_p,t}(s_{u_p})=e_{u_p,t}(s_{u_p})\,\partial_{u_j}s_{u_p}\int_{u_p}^t\partial_x^2 g(r,s_r)e_{u_p,r}(s_{u_p})\,dr.$$

*Then for* $p \geq j$ *and* $u_p \leq t < u_{p+1}$*, we have the recurrence relations*

$$\begin{aligned}\partial_{u_j}^2 s_t &= e_{u_p,t}(s_{u_p})\,\partial_{u_j}^2 s_{u_p} + \rho_{j,p}(t,u,a)\,\partial_{u_j} s_{u_p},\\ \partial_{u_j}^2 s_{u_{p+1}} &= \partial_x^2 c(u_{p+1},a_{p+1},s_{u_{p+1}-})(e_{u_p,u_{p+1}}(s_{u_p})\,\partial_{u_j} s_{u_p})^2\\ &\quad + (1+\partial_x c)(u_{p+1},a_{p+1},s_{u_{p+1}-})\\ &\qquad \times (\rho_{j,p}(u_{p+1})\,\partial_{u_j} s_{u_p} + e_{u_p,u_{p+1}}(s_{u_p})\,\partial_{u_j}^2 s_{u_p}).\end{aligned}$$

B. Derivatives with respect to $a_j$. *For* $t < u_j, \partial_{a_j} s_{u_j}(u,a) = 0$ *and for* $t \geq u_j$, $\partial_{a_j} s_t(u,a)$ *satisfies the equation*

$$\begin{aligned}\partial_{a_j} s_t &= \partial_a c(u_j,a_j,s_{u_j-}) + \sum_{i=j+1}^{J_t(u)} \partial_x c(u_i,a_i,s_{u_i-})\,\partial_{a_j} s_{u_i-}\\ &\quad + \int_{u_j}^t \partial_x g(r,s_r)\,\partial_{a_j} s_r\,dr.\end{aligned}\tag{3.7}$$

*The second-order derivatives satisfy*

$$\begin{aligned}\partial_{a_j}^2 s_t &= \partial_a^2 c(u_j,a_j,s_{u_j-}) + \sum_{i=j+1}^{J_t(u)} \partial_x^2 c(u_i,a_i,s_{u_i-})(\partial_{a_j} s_{u_i-})^2\\ &\quad + \int_{u_j}^t \partial_x^2 g(r,s_r)(\partial_{a_j} s_r)^2\,dr\\ &\quad + \sum_{i=j+1}^{J_t(u)} \partial_x c(u_i,a_i,s_{u_i-})\,\partial_{a_j}^2 s_{u_i-} + \int_{u_j}^t \partial_x g(r,s_r)\,\partial_{a_j}^2 s_r\,dr.\end{aligned}\tag{3.8}$$

Proof. A. It is clear that for $t < u_j$, $s_t$ does not depend on $u_j$ and so $\partial_{u_j} s_t = 0$. We now compute

$$\partial_{u_j} s_{u_j-} = \partial_{u_j}\Phi_{u_{j-1}}(u_j,s_{u_{j-1}}) = g(u_j,\Phi_{u_{j-1}}(u_j,s_{u_{j-1}})) = g(u_j,s_{u_j-}).$$

Then

$$\begin{aligned}\partial_{u_j} s_{u_j} &= \partial_{u_j}(s_{u_j-} + c(u_j,a_j,s_{u_j-}))\\ &= \partial_t c(u_j,a_j,s_{u_j-}) + (1+\partial_x c(u_j,a_j,s_{u_j-}))\,\partial_{u_j} s_{u_j-}\\ &= \partial_t c(u_j,a_j,s_{u_j-}) + (1+\partial_x c(u_j,a_j,s_{u_j-}))g(u_j,s_{u_j-}).\end{aligned}$$

For $u_j < t < u_{j+1}$, we have

$$\begin{aligned}\partial_{u_j} s_t &= \partial_{u_j}\Phi_{u_j}(t,s_{u_j}) = e_{u_j,t}(s_{u_j})(-g(u_j,s_{u_j}) + \partial_{u_j} s_{u_j})\\ &= e_{u_j,t}(s_{u_j})(-g(u_j,s_{u_j}) + \partial_t c(u_j,a_j,s_{u_j-})\\ &\qquad + (1+\partial_x c(u_j,a_j,s_{u_j-}))g(u_j,s_{u_j-}))\\ &= e_{u_j,t}(s_{u_j})q(u_j,a_j,s_{u_j-})\end{aligned}$$

and the same computation gives $\partial_{u_j} s_{u_{j+1}-} = e_{u_j,u_{j+1}}(s_{u_j}) q(u_j, a_j, s_{u_j-})$. Finally,

$$\begin{aligned}\partial_{u_j} s_{u_{j+1}} &= (1 + \partial_x c(u_{j+1}, a_{j+1}, s_{u_{j+1}-}))\, \partial_{u_j} s_{u_{j+1}-}\\ &= (1 + \partial_x c(u_{j+1}, a_{j+1}, s_{u_{j+1}-})) e_{u_j,u_{j+1}}(s_{u_j}) q(u_j, a_j, s_{u_j-}).\end{aligned}$$

We now assume that $u_p \leq t < u_{p+1}$, $p \geq j + 1$, and we write

$$\partial_{u_j} s_t = \partial_{u_j} \Phi_{u_p}(t, s_{u_p}) = e_{u_p,t}(s_{u_p})\, \partial_{u_j} s_{u_p}.$$

The same computation gives $\partial_{u_j} s_{u_{p+1}-} = e_{u_p,u_{p+1}}(s_{u_p})\, \partial_{u_j} s_{u_p}$. Finally, we have

$$\begin{aligned}\partial_{u_j} s_{u_p} &= \partial_{u_j}(s_{u_p-} + c(u_p, a_p, s_{u_p-})) = (1 + \partial_x c(u_p, a_p, s_{u_p-}))\, \partial_{u_j} s_{u_p-}\\ &= (1 + \partial_x c(u_p, a_p, s_{u_p-})) e_{u_{p-1},u_p}(s_{u_{p-1}})\, \partial_{u_j} s_{u_{p-1}}\end{aligned}$$

and the proof is complete for the first-order derivatives. The relations concerning the second-order derivatives are obtained by direct computations.

B. Using the recurrence relations (3.4), one verifies that for every $t \in [0, T]$, $a_j \to s_t(u, a)$ is continuously differentiable and one may then differentiate in equation (3.2) (this was not possible in the case of the derivatives with respect to $u_j$ because these derivatives are not continuous). □

As an immediate consequence of the above lemma we obtain:

COROLLARY 3.1. *Suppose that Hypothesis 3.1 holds true and that the starting point $x$ satisfies $|x| \leq K$ for some $K$. Then for each $n \in \mathbb{N}$ and $T > 0$, there exists a constant $C_n(K,T)$ such that for every $0 < u_1 < \cdots < u_n < T$, $a \in \mathbb{R}^n$ and $0 \leq t \leq T$,*

$$\max_{j=1,\ldots,n} (|s_t| + |\partial_{u_j} s_t| + |\partial^2_{u_j} s_t| + |\partial_{a_j} s_t| + |\partial^2_{a_j} s_t|)(u, a) \leq C_n(K, T). \tag{3.9}$$

Finally, we present a corollary which is useful in order to control the nondegeneracy.

COROLLARY 3.2. *Assume that Hypothesis 3.1 holds true and there exists a constant $\eta > 0$ such that for every $(t, a, x) \in [0, T] \times \mathbb{R} \times \mathbb{R}$, one has*

$$\begin{aligned}|1 + \partial_x c(t, a, x)| &\geq \eta,\\ |q(t, a, x)| &\geq \eta.\end{aligned} \tag{3.10}$$

*Let $n \in \mathbb{N}$ be fixed. Then there exists a constant $\varepsilon_n > 0$ such that for every $j = 1, \ldots, n$ and every $(u, a) \in [0, T]^n \times \mathbb{R}^n$, we have,*

$$\inf_{t > u_j} |\partial_{u_j} s_t(u, a)| \geq \varepsilon_n. \tag{3.11}$$

PROOF. Since $\partial_x g$ is bounded, there exists a constant $C > 0$ such that $e_{s,t}(x) \geq e^{-CT}$ for $0 \leq s < t \leq T$. One then employs (3.5) and (3.6). □

3.2. *Integration by parts with respect to the amplitudes of the jumps.* In this section, we will use the integration by parts formula for $S_t$ which will be regarded as a simple functional of $\Delta_i, i\in\mathbb{N}$. So, with the notation from Section 2, we have $V_i=\Delta_i$ and $\mathcal{G}=\sigma\{T_i : i\in\mathbb{N}\}$. We recall that $J_t=n$ on $\{T_n\le t<T_{n+1}\}$. Then, on $\{J_t=n\}$, we have

$$S_t = s_t(T_1,\dots,T_n,\Delta_1,\dots,\Delta_n),$$

where $s_t$ is defined in the previous section [see (3.2)].

We assume that Hypotheses 3.1 and 2.1 [i.e., $\mathrm{E}(|\Delta_i|^p)<\infty$ for all $p\in\mathbb{N}$] hold true. Moreover, we consider some $q_0<q_1<\cdots<q_{k+1}$ and define

$$I=\bigcup_{i=0}^{k}(q_i,q_{i+1}).$$

We assume the following.

Hypothesis 3.2. The law of $\Delta_i$ is absolutely continuous on $I$ with respect to the Lebesgue measure and has the density $p(y)=\mathbf{1}_I(y)e^{\rho(y)}$, that is,

$$\mathrm{E}(f(\Delta_i))=\int_I f(y)e^{\rho(y)}\,dy$$

for every measurable and positive function $f$.

The function $\rho$ is assumed to be continuously differentiable and bounded on $I$.

Therefore, Hypothesis 2.2 holds true.

Since $\rho$ is not differentiable on the whole of $\mathbb{R}$, we work with the following weight. We take $\alpha\in(0,1)$ and $\beta>\alpha$ and we define

$$\pi(y)=\begin{cases}(q_{i+1}-y)^{\alpha}(y-q_i)^{\alpha}, & \text{for } y\in(q_i,q_{i+1}), i=0,\dots,k,\\ 0, & \text{for } y\in(q_0,q_{k+1})^c.\end{cases}$$

We introduce the following convention: if $b=q_{k+1}=+\infty$ or $a=q_0=-\infty$, we define

$$\pi(y)=\begin{cases}(y-q_k)^{\alpha}|y|^{-\beta}, & \text{for } y>q_k,\\ (q_1-y)^{\alpha}|y|^{-\beta}, & \text{for } y<q_1.\end{cases}$$

Since $\alpha\in(0,1)$, we can show by elementary computations that $\pi$ satisfies Hypothesis 2.3, that is, $\pi\in L_{(\infty)}(A)$ and $\pi'\in L_{(1+)}(A)$.

Let $A:=\{J_t=n\}$. In view of (3.9),

$$(a_1,\dots,a_n)\to s_t(T_1(\omega),\dots,T_n(\omega),a_1,\dots,a_n)$$

is twice continuously differentiable and has bounded derivatives. It follows that $S_t\in S_{(n,2)}(A)$.

The differential operators which appear in the integration by parts formula are

$$D_i S_t = \partial_{a_i} s_t(T_1,\dots,T_n,\Delta_1,\dots,\Delta_n),$$

$$LS_t = -\sum_{i=1}^{n}\biggl(\pi(\Delta_i)\partial_{a_i}^2 s_t(T_1,\dots,T_n,\Delta_1,\dots,\Delta_n) + \Bigl(\pi' + \pi\frac{\rho'}{\rho}\Bigr)(\Delta_i)\,\partial_{a_i} s_t(T_1,\dots,T_n,\Delta_1,\dots,\Delta_n)\biggr),$$

$$\sigma_{S_t} = \sum_{i=1}^{n}\pi(\Delta_i)|D_i S_t|^2 = \sum_{i=1}^{n}\pi(\Delta_i)|\partial_{a_i} s_t(T_1,\dots,T_n,\Delta_1,\dots,\Delta_n)|^2,$$

$$\gamma_{S_t} = \frac{1}{\sigma_{S_t}} = \frac{1}{\sum_{i=1}^{n}\pi(\Delta_i)|\partial_{a_i} s_t(T_1,\dots,T_n,\Delta_1,\dots,\Delta_n)|^2}.$$

All of these quantities may be computed using (3.7) and (3.8).

The result which is used in sensitivity computations is the following.

PROPOSITION 3.1. *Suppose that Hypotheses 3.1 and 3.2 hold true and, moreover, suppose that there exists a positive constant $\eta$ such that for every $t,a,x$, we have*

$$\text{(3.12)}\qquad \begin{aligned} &\text{(i)} && |\partial_a c(t,a,x)| \ge \eta,\\ &\text{(ii)} && |1+\partial_x c(t,a,x)| \ge \eta. \end{aligned}$$

*Take $\alpha \in (0,1/2)$ and $\beta > \alpha$. Then for every differentiable function $\phi:\mathbb{R}\to\mathbb{R}$ which has linear growth and for every $n \ge 1$,*

$$\mathrm{E}(\phi'(S_t)\,\partial_x S_t \mathbf{1}_{\{J_t=n\}}) = \mathrm{E}(\phi(S_t) H_n \mathbf{1}_{\{J_t=n\}}),$$

*with, on $\{J_t = n\}$,*

$$\begin{aligned} H_n &:= H_n(S_t,\partial_x S_t)\\ &= \partial_x S_t \gamma_{S_t} LS_t - \gamma_{S_t}\langle DS_t, D(\partial_x S_t)\rangle_\pi - \partial_x S_t \langle DS_t, D\gamma_{S_t}\rangle_\pi. \end{aligned}$$

PROOF. Let $n \in \mathbb{N}^*$ be fixed. We already know that $S_t \in \mathcal{S}_{(n,2)}(A)$ with $A = \{J_t = n\}$.

Moreover, $\partial_x S_t = \partial_x s_t(T_1,\dots,T_n,\Delta_1,\dots,\Delta_n)$ and $\partial_x s_t$ is computed by the recurrence relations

$$\partial_x s_0 = 1,$$

$$\partial_x s_t = (1+\partial_x c(u_i,a_i,s_{u_i-}))\,\partial_x s_{u_i-} + \int_{u_i}^{t}\partial_x g(r,s_r)\,\partial_x s_r\,dr, \qquad u_i \le t < u_{i+1}.$$

It is then easy to check that $\partial_x s_t$ and its derivatives with respect to $a_i$, $i=1,\dots,n$, are bounded on the set $\{J_t=n\}$ and consequently that $\partial_x S_t \in \mathcal{S}_{(n,1)}(A)$.

• Suppose that $n=1$. We will use Corollary 2.2, so we check that the nondegeneracy condition (2.12) holds true. One has

$$\partial_{a_1} s_t = \partial_a c(u_1, a_1, s_{u_1-}) + \int_{u_1}^t \partial_x g(r, s_r)\, \partial_{a_1} s_r\, dr$$

so that, using (3.12),

$$|\partial_{a_1} s_t| = |\partial_a c(u_1, a_1, s_{u_1-})| \exp\left(\int_{u_1}^t \partial_x g(r, s_r)\right) \geq c$$

for some positive constant $c$. Inequality (2.12) then follows. Then the integration by parts formula (2.13) holds true for $S_t$ and $\partial_x S_t$ on the event $A=\{J_t=1\}$. Moreover, by our choice of $\pi$, the border terms are canceled, which gives $\mathrm{E}(\phi'(S_t)\,\partial_x S_t \mathbf{1}_{\{J_t=1\}}) = \mathrm{E}(\phi(S_t) H_1 \mathbf{1}_{\{J_t=1\}})$ with

$$\begin{aligned} H_1 \mathbf{1}_{\{J_t=1\}} &= \delta(\partial_x S_t \gamma_t D S_t) \mathbf{1}_{\{J_t=1\}} \\ &= -\partial_{a_1}(\pi(\Delta_1)\,\partial_x S_t \gamma_{S_t} D S_t) - \pi(\Delta_1)\,\partial \ln p \gamma_{S_t} D S_t\, \partial_x S_t \mathbf{1}_{\{J_t=1\}}. \end{aligned}$$

On $\{J_t=1\}$, we have

$$\begin{aligned} \pi(\Delta_1)\,\partial_x S_t \gamma_{S_t} D S_t &= \frac{\pi(\Delta_1)\,\partial_x s_t(T_1,\Delta_1)\,\partial_{a_1} s_t(T_1,\Delta_1)}{\pi(\Delta_1)|\partial_{a_1} s_t(T_1,\Delta_1)|^2} \\ &= \frac{\partial_x s_t(T_1,\Delta_1)}{\partial_{a_1} s_t(T_1,\Delta_1)} \mathbf{1}_I(\Delta_i). \end{aligned}$$

• Now suppose that $n \geq 2$. In this case, we will use Theorem 2.2, so we look at the nondegeneracy condition (2.11). Since $\pi$ is bounded, this amounts to finding $\delta > 0$ such that for $i=1,\dots,n$,

$$\mathrm{E}[\mathbf{1}_{\{J_t=n\}}((1+|\pi'(\Delta_i)|)\gamma_{S_t}^2)^{1+\delta}] < \infty. \tag{3.13}$$

We recall that $\pi(y) = \sum_{i=0}^k (q_{i+1}-y)^\alpha (y-q_i)^\alpha \mathbf{1}_{(q_i,q_{i+1})}(y)$, so

$$\pi'(y) = \begin{cases} \alpha(q_{i+1}-y)^{\alpha-1}(y-q_i)^{\alpha-1}(q_i - 2y + q_{i+1}), & \text{if } y \in (q_i, q_{i+1}), \\ 0, & \text{if } y \in (q_0, q_{k+1})^c. \end{cases}$$

We choose $\delta > 0$ such that $2\alpha(1+\delta) < 1$ and $(1-\alpha)(1+\delta) < 1$ [which is possible because $\alpha \in (0, 1/2)$]. In particular, since $\rho$ is bounded on $I$ and $\Delta_i$ have finite moments of any order, this gives

$$\mathrm{E}(\pi(\Delta_i)^{-2(1+\delta)}) < \infty \quad \text{and} \quad \mathrm{E}(|\pi'(\Delta_i)|^{1+\delta}) < \infty.$$

The proof of (3.13) is different for $i=n$ and $i=1,\dots,n-1$.

First, take $i < n$. One has $|\partial_{a_n} s_t| = |\partial_a c(u_n, a_n, s_{u_n-})| \exp(\int_{u_n}^t \partial_x g(r, s_r)) \geq c$, so $\sigma_{S_t} \geq c^2\pi^2(\Delta_n)$. Since $\Delta_i$ and $\Delta_n$ are independent,

$$\begin{aligned}\mathrm{E}[\mathbf{1}_{\{J_t=n\}}((1+|\pi'(\Delta_i)|)\gamma_{S_t}^2)^{1+\delta}] &\leq c^{-2}\mathrm{E}[\mathbf{1}_{\{J_t=n\}}((1+|\pi'(\Delta_i)|)\pi^{-2}(\Delta_n))^{1+\delta}] \\ &= c^{-2}\mathrm{E}(\pi^{-2(1+\delta)}(\Delta_n))\mathrm{E}[(1+|\pi'(\Delta_i)|)^{1+\delta}] \\ &< \infty.\end{aligned}$$

Now, take $i = n$ and write $\sigma_{S_t} \geq \pi^2(\Delta_{n-1})|D_{n-1}S_t|^2$. A simple computation shows that

$$\partial_{a_{n-1}} s_t = \partial_a c(u_{n-1}, a_{n-1}, s_{u_{n-1}^-})(1+\partial_x c(u_n, a_n, s_{u_n^-}))\exp\left(\int_{u_{n-1}}^t \partial_x g(r, s_r)\,dr\right).$$

Using (3.12), we obtain $\partial_{a_{n-1}} s_t \geq c > 0$, so $\sigma_{S_t} \geq c^2\pi^2(\Delta_{n-1})$. Consequently,

$$\begin{aligned}&\mathrm{E}[\mathbf{1}_{\{J_t=n\}}((1+|\pi'(\Delta_n)|)\gamma_{S_t}^2)^{1+\delta}] \\ &\quad\leq c^{-2}\mathrm{E}[\mathbf{1}_{\{J_t=n\}}((1+|\pi'(\Delta_i)|)\pi^{-2}(\Delta_n))^{1+\delta}] \\ &\quad= c^{-2}\mathrm{E}(\pi^{-2(1+\delta)}(\Delta_{n-1}))\mathrm{E}[(1+|\pi'(\Delta_n)|)^{1+\delta}] < \infty\end{aligned}$$

and the proof is complete. □

REMARK 3.1. Suppose now that $\rho$ is differentiable on the whole of $\mathbb{R}$. We then take no weight, $\pi = 1$ and hypothesis (3.12)(i) gives $\sigma_{S_t} \geq c$ on $\{J_t = n\}$ for all $n \in \mathbb{N}^*$. So the above integrability problems disappear. In particular, hypothesis (3.12)(ii) is no longer necessary. This case is discussed in [1].

3.3. *Integration by parts with respect to the jump times.* In this section, we differentiate with respect to the jump times $T_i$, $i \in \mathbb{N}$. It is well known (see [2]) that conditionally to $\{J_t = n\}$, the law of the vector $(T_1, \ldots, T_n)$ is absolutely continuous with respect to the Lebesgue measure and has the following density:

$$p(\omega, t_1, \ldots, t_n) = \frac{n!}{t^n}\mathbf{1}_{\{0<t_1<\cdots<t_n<t\}}(t_1, \ldots, t_n)\mathbf{1}_{\{J_t(\omega)=n\}}.$$

In particular, for a given $i = 1, \ldots, n$, conditionally to $\{J_t = n\}$ and to $\{T_j, j \neq i\}$, $T_i$ is uniformly distributed on $[T_{i-1}(\omega), T_{i+1}(\omega)]$. Therefore, it has the density (with the convention $T_0 = 0, T_{n+1} = t$)

$$p_i(\omega, u) = \frac{1}{T_{i+1}(\omega) - T_{i-1}(\omega)}\mathbf{1}_{[T_{i-1}(\omega), T_{i+1}(\omega)]}(u)\,du, \qquad i = 1, \ldots, n.$$

Since $p_i$ is not differentiable with respect to $u$, we must use the following weights:

$$\pi_i(\omega, u) = (T_{i+1}(\omega) - u)^\alpha(u - T_{i-1}(\omega))^\alpha\mathbf{1}_{[T_{i-1}(\omega), T_{i+1}(\omega)]}(u), \qquad i = 1, \ldots, n,$$

with $\alpha \in (0, 1/2)$.

In order to fit with the notation from the first section, we take $V_i = T_i$, $k_i = 2$, $t_i^1 = T_{i-1}$ and $t_i^2 = T_{i+1}$. We have $\mathcal{G} = \sigma(\Delta_i, i \in \mathbb{N}) \vee \sigma(J_t)$. We fix $n$ and work on the set $A = \{J_t = n\}$. Hypotheses 2.1, 2.2 and 2.3 then hold true and $S_t = s_t(T_1, \ldots, T_n, \Delta_1(\omega), \ldots, \Delta_n(\omega))$. So $S_t$ is a simple functional and the function which represents $S_t$ is twice differentiable and has continuous derivatives on the whole of $\mathbb{R}^n$. The differential operators are

$$D_i S_t = \partial_{u_i} s_t(T_1, \ldots, T_n, \Delta_1(\omega), \ldots, \Delta_n(\omega)),$$

$$\sigma_{S_t} = \sum_{i=1}^n \pi_i(\omega, T_i) |\partial_{u_i} s_t(T_1, \ldots, T_n, \Delta_1(\omega), \ldots, \Delta_n(\omega))|^2,$$

$$L_i S_t = -(\pi_i' \partial_{u_i} s_t + \pi_i \partial_{u_i}^2 s_t)(T_1, \ldots, T_n, \Delta_1(\omega), \ldots, \Delta_n(\omega))$$

and all of these quantities may be computed using Lemma 3.1.

PROPOSITION 3.2. *Suppose that Hypothesis 3.1 holds true. Suppose, moreover, that (3.10) is satisfied, that is, that*

$$|q(t, a, x)| \geq \eta > 0,$$
$$|(1 + \partial_x c)(t, a, x)| \geq \eta > 0$$

*for some $\eta > 0$. Take $\alpha \in (0, \frac{1}{2})$. Then for every $n \geq 4$ and every continuously differentiable function $\phi$ which has linear growth, we have*

$$\mathrm{E}(\phi'(S_t)\, \partial_x S_t \mathbf{1}_{\{J_t = n\}}) = \mathrm{E}(\phi(S_t) H_n \mathbf{1}_{\{J_t = n\}})$$

*with, on $\{J_t = n\}$,*

$$\begin{aligned} H_n &:= H_n(S_t, \partial_x S_t) \\ &= \partial_x S_t \gamma_{S_t} L S_t - \gamma_{S_t} \langle DS_t, D(\partial_x S_t)\rangle_\pi - \partial_x S_t \langle DS_t, D\gamma_{S_t}\rangle_\pi. \end{aligned}$$

PROOF. From (3.9), we know that $s_t(u, a)$ and its derivatives up to order two with respect to $u_i$, $i = 1, \ldots, n$, are bounded on $[0, T]^n$. It follows that $S_t \in \mathcal{S}_{(n,2)}(A)$.

Since $\pi_i$ are bounded, the nondegeneracy condition (2.11) amounts to

$$\mathrm{E}[\mathbf{1}_{\{J_t = n\}} \gamma_{S_t}^{2(1+\eta)}] < \infty \quad \text{and} \quad \mathrm{E}[\mathbf{1}_{\{J_t = n\}} \gamma_{S_t}^{2(1+\eta)} |\pi_i'(T_i)|^{1+\eta}] < \infty$$

for some $\eta > 0$.

Let us prove that $\mathrm{E}[\mathbf{1}_{\{J_t = n\}} \gamma_{S_t}^{2(1+\eta)} |\pi_i'(T_i)|^{1+\eta}] < \infty$. We define $\delta_i = T_i - T_{i-1}$ and $\delta_{n+1} = T - T_n$, so that $\pi_i = \delta_i^\alpha \delta_{i+1}^\alpha$. We use (3.11) in order to obtain

$$\sigma_{S_t} = \sum_{i=1}^n \delta_{i+1}^\alpha \delta_i^\alpha |\partial_{u_i} s_t(T_1, \ldots, T_n, \Delta_1, \ldots, \Delta_n)|^2 \geq \varepsilon^2 \sum_{i=1}^n \delta_{i+1}^\alpha \delta_i^\alpha.$$

Since $\pi_i'(T_i) = \alpha(-\delta_{i+1}^{\alpha-1}\delta_i^{\alpha} + \delta_{i+1}^{\alpha}\delta_i^{\alpha-1})$, we must check that, for every $i = 1, \dots, n$,

$$\mathrm{E}\left[(\delta_i^{\alpha-1}\delta_{i+1}^{\alpha} + \delta_i^{\alpha}\delta_{i+1}^{\alpha-1})^{1+\eta}\left(\sum_{j=1}^{n}\delta_{j+1}^{\alpha}\delta_j^{\alpha}\right)^{-2(1+\eta)}\right] < \infty.$$

Take $i = 1$ and write

$$\begin{aligned}\mathrm{E}\left[(\delta_1^{\alpha-1}\delta_2^{\alpha})^{1+\eta}\left(\sum_{j=1}^{n}\delta_{j+1}^{\alpha}\delta_j^{\alpha}\right)^{-2(1+\eta)}\right] &\le \mathrm{E}[(\delta_1^{\alpha-1}\delta_2^{\alpha})^{1+\eta}(\delta_2^{\alpha}\delta_3^{\alpha})^{-2(1+\eta)}] \\ &= \mathrm{E}(\delta_1^{(\alpha-1)(1+\eta)})\mathrm{E}(\delta_2^{-\alpha(1+\eta)})\mathrm{E}(\delta_3^{-2\alpha(1+\eta)}).\end{aligned}$$

Since $\delta_i$ is exponentially distributed of parameter $\mu(\mathbb{R})$, a necessary and sufficient condition in order to have $\mathrm{E}(\delta_i^{-p}) < \infty$ is $p < 1$. We then choose $\eta$ sufficiently small that $(1-\alpha)(1+\eta) < 1$ and $\alpha(1+\eta) < 2\alpha(1+\eta) < 1$ (which is possible because $0 < \alpha < 1/2$) and we have $\mathrm{E}(\delta_1^{(\alpha-1)(1+\eta)}) < \infty$, $\mathrm{E}(\delta_2^{-\alpha(1+\eta)}) < \infty$ and $\mathrm{E}(\delta_3^{-2\alpha(1+\eta)}) < \infty$. So

$$\mathrm{E}\left[(\delta_1^{\alpha-1}\delta_2^{\alpha})^{1+\eta}\left(\sum_{j=1}^{n}\delta_{j+1}^{\alpha}\delta_j^{\alpha}\right)^{-2(1+\eta)}\right] < \infty.$$

We now write

$$\begin{aligned}&\mathrm{E}\left[(\delta_1^{\alpha}\delta_2^{\alpha-1})^{1+\eta}\left(\sum_{j=1}^{n}\delta_{j+1}^{\alpha}\delta_j^{\alpha}\right)^{-2(1+\eta)}\right] \\ &\quad\le \mathrm{E}[(\delta_1^{\alpha}\delta_2^{\alpha-1})^{1+\eta}(\delta_3^{\alpha}\delta_4^{\alpha})^{-2(1+\eta)}] \\ &\quad= \mathrm{E}(\delta_2^{(\alpha-1)(1+\eta)})\mathrm{E}(\delta_1^{\alpha(1+\eta)})\mathrm{E}(\delta_3^{-2\alpha(1+\eta)})\mathrm{E}(\delta_4^{-2\alpha(1+\eta)}).\end{aligned}$$

Recalling that $\delta_i$ has finite moments of any order, by the choice of $\eta$, we obtain

$$\mathrm{E}\left[(\delta_1^{\alpha}\delta_2^{\alpha-1})^{1+\eta}\left(\sum_{j=1}^{n}\delta_{j+1}^{\alpha}\delta_j^{\alpha}\right)^{-2(1+\eta)}\right] < \infty.$$

Since $n \ge 4$, the same argument works for $i = 2, \dots, n$ and leads to $\mathrm{E}[\mathbf{1}_{\{J_t = n\}} \times \gamma_{S_t}^{2(1+\eta)}] < \infty$. □

Remark 3.2. Suppose that $n = 2$. Then

$$\begin{aligned}&(\delta_1^{\alpha-1}\delta_2^{\alpha})^{1+\eta}\left(\sum_{j=1}^{n}\delta_{j+1}^{\alpha}\delta_j^{\alpha}\right)^{-2(1+\eta)} \\ &\quad= (\delta_1^{\alpha-1}\delta_2^{\alpha})^{1+\eta}\delta_2^{-2\alpha(1+\eta)}(\delta_1^{\alpha} + \delta_3^{\alpha})^{-2(1+\eta)} \\ &\quad= \delta_2^{-\alpha(1+\eta)} \times (\delta_1^{-(\alpha+1)(1+\eta)} + \delta_3^{-2\alpha(1+\eta)}\delta_1^{-(1-\alpha)(1+\eta)})\end{aligned}$$

and this quantity is not integrable for $\alpha > 0$, $\eta > 0$.

REMARK 3.3. For $n = 1$, one may use Corollary 2.2 in order to obtain an integration by parts formula.

But for $n = 2$ and $n = 3$, we are not able to handle the nondegeneracy problem. In our numerical examples, we will use the noise coming from the amplitudes of the jumps in order to solve the problem for $n = 2$ and $n = 3$.

3.3.1. *Examples.* • We consider the geometrical model

$$dS_t = S_t(r dt + \alpha(t,a)\, dN(t,a)).$$

In this case, $g(t,x) = xr$ and $c(t,a,x) = x\alpha(t,a)$. It follows that

$$q(t,a,x) = x\,\partial_t\alpha(t,a) + xr\alpha(t,a) + xr - r(x + x\alpha(t,a)) = x\,\partial_t\alpha(t,a).$$

In particular, if $\alpha$ does not depend on the time, the model is degenerate from the point of view of the jump times. The nondegeneracy condition becomes

$$|\partial_t\alpha(t,a)| \geq \varepsilon.$$

On the other hand, the condition $|1 + \partial_x c(t,a,x)| \geq \eta$ becomes

$$|1 + \alpha(t,a)| \geq \eta.$$

• We now consider the following Vasicek-type model:

$$dS_t = S_t r\, dt + \alpha(t,a)\, dN(t,a).$$

In this case, $g(t,x) = xr$ and $c(t,a,x) = \alpha(t,a)$. It follows that

$$q(t,a,x) = \partial_t\alpha(t,a) + xr - r(x + \alpha(t,a)) = \partial_t\alpha(t,a) - r\alpha(t,a).$$

Suppose that $\alpha$ does not depend on the time so that $\partial_t\alpha = 0$. Then the nondegeneracy assumption becomes

$$|\alpha(a)| \geq \varepsilon.$$

The condition $|1 + \partial_x c(t,a,x)| \geq \eta$ becomes

$$|1 + \alpha(a)| \geq \eta.$$

3.4. *Mixed calculus.* In this section, we briefly present the differential calculus with respect to both noises coming from the jump amplitudes and from the jump times. So the random variables will be $V_i = T_i, i = 1,\dots,n$, $V_{n+i} = \Delta_i$, $i = 1,\dots,n$, and $\mathcal{G} = \sigma(J_t)$. We combine the results from the two previous sections (and we keep the notation therefrom). We still assume Hypotheses 3.1 and 3.2. The differential operators are

$$D_i S_t = \begin{cases} \partial_{u_i} s_t(u_1,\dots,u_n,\Delta_1(\omega),\dots,\Delta_n(\omega)), & i = 1,\dots,n, \\ \partial_{a_{i-n}} s_t(T_1,\dots,T_n,\Delta_1,\dots,\Delta_n), & i = n+1,\dots,2n. \end{cases}$$

We will use the weights defined in the previous sections, namely

$$\pi_i(\omega,u)=(T_{i+1}(\omega)-u)^\alpha(u-T_{i-1}(\omega))^\alpha\mathbf{1}_{[T_{i-1}(\omega),T_{i+1}(\omega)]}(u),\qquad i=1,\dots,n,$$

$$\pi_i(y)=\pi(y)=\sum_{p=1}^{k-1}(q_{p+1}-y)^\alpha(y-q_p)^\alpha\mathbf{1}_{(q_p,q_{p+1})}(y),\qquad i=n+1,\dots,2n,$$

where $\alpha\in(0,\frac{1}{2})$.

We have

$$L_iS_t=\begin{cases}-(\pi_i'(T_i)\,\partial_{u_i}s_t+\pi_i(T_i)\,\partial_{u_i}^2s_t), & \text{for } i=1,\dots,n,\\ -(\pi(\Delta_i)\,\partial_{a_i}^2s_t+(\pi'+\pi\rho')(\Delta_i)\,\partial_{a_i}s_t), & \text{for } i=n+1,\dots,2n.\end{cases}$$

Finally, $LS_t=\sum_{i=1}^{2n}L_iS_t$. All of these quantities may be computed using the formulas from the previous sections.

THEOREM 3.1. *Suppose that Hypotheses 3.1 and 3.2 hold true and that*

$$\text{(i)}\qquad |q(t,a,x)|\ge\varepsilon>0,$$

$$\text{(ii)}\qquad |\partial_a c(t,a,x)|\ge\varepsilon>0,$$

$$\text{(iii)}\qquad |(1+\partial_x c)(t,a,x)|\ge\varepsilon>0.$$

*Then for every $n\ge1$ and every continuously differentiable function $\phi$ which has linear growth, we have*

$$\mathrm{E}(\phi'(S_t)\,\partial_xS_t\mathbf{1}_{\{J_t=n\}})=\mathrm{E}(\phi(S_t)H_n\mathbf{1}_{\{J_t=n\}})$$

*with, on $\{J_t=n\}$,*

$$\begin{aligned}H_n&:=H_n(S_t,\partial_xS_t)\\&=\partial_xS_t\gamma_{S_t}LS_t-\gamma_{S_t}\langle DS_t,D(\partial_xS_t)\rangle_\pi-\partial_xS_t\langle DS_t,D\gamma_{S_t}\rangle_\pi.\end{aligned}$$

PROOF. We write

$$\begin{aligned}\sigma_{S_t}&\ge(\pi_n(\omega,T_n)|\partial_{u_n}s_t|^2+\pi(\Delta_n)|\partial_{a_n}s_t|^2)(T_1,\dots,T_n,\Delta_1,\dots,\Delta_n)\\&\ge\varepsilon^2(\pi_n(\omega,T_n)+\pi(\Delta_n))\end{aligned}$$

for some $\varepsilon>0$. Then, using the same techniques as in the proofs of Propositions 3.2 and 3.1, one shows that the nondegeneracy (2.11) condition holds true. □

REMARK 3.4. Note that the nondegeneracy condition holds true for every $n$ (including $n=2$) because we may use the noises coming from the jump times and the jump amplitudes at the same time.

## 4. Numerical results.

4.1. *Malliavin estimator.* In this section, we compute the Delta of two European options: call option with payoff $\phi(x)=(x-K)_+$ and digital option with payoff $\phi(x)=\mathbf{1}_{x\geq K}$. The asset $(S_t)_{t\geq 0}$ follows a one-dimensional pure jump diffusion process. We use the notation from the beginning of Section 3.

We deal with two different pure jump diffusion models. The first is a Vasicek-type model,

$$S_t = x - \int_0^t r(S_u-\alpha)\,du + \sum_{i=1}^{J_t}\sigma\Delta_i, \tag{4.1}$$

and the second is a geometrical model,

$$S_t = x + \int_0^t rS_u\,du + \sigma\sum_{i=1}^{J_t} S_{T_i^-}\Delta_i. \tag{4.2}$$

In both models, we take $\Delta_i\sim\mathcal{N}(0,1)$, $i\geq 1$. That is, for all $i\geq 1$, $\Delta_i$ has the density $p(x)=\frac{1}{\sqrt{2\pi}}e^{\rho(x)}$ with $\rho(x)=-\frac{x^2}{2}$. Note that even if $\rho$ is not bounded on $\mathbb{R}$, the integration by parts formula holds by a truncation argument.

Our aim is to compute $\partial_x\mathrm{E}(\phi(S_T))$ using the integration by parts formula derived in the previous sections. We write

$$\begin{aligned}\partial_x\mathrm{E}(\phi(S_T)) &= \mathrm{E}(\phi'(S_T)\,\partial_x S_T)\\ &= \mathrm{E}(\phi'(S_T)\,\partial_x S_T\mathbf{1}_{\{J_T=0\}}) + \sum_{n=1}^{\infty}\mathrm{E}(\phi'(S_T)\,\partial_x S_T\mathbf{1}_{\{J_T=n\}}).\end{aligned}$$

For $n\geq 1$, we use the integration by parts formula on $\{J_T=n\}$ and obtain

$$\mathrm{E}(\phi'(S_T)\,\partial_x S_T\mathbf{1}_{\{J_T=n\}}) = \mathrm{E}(\phi(S_T)H_n\mathbf{1}_{\{J_T=n\}}),$$

where $H_n$ is a weight involving Malliavin derivatives of $S_T$ and $\partial_x S_T$. Summing over $n=1,2,\dots,$ we obtain

$$\partial_x\mathrm{E}(\phi(S_T)) = \mathrm{E}(\phi'(S_T)\partial_x S_T\mathbf{1}_{\{J_T=0\}}) + \mathrm{E}(\phi(S_T)H_{J_T}(S_T,\partial_x S_T)\mathbf{1}_{\{J_T\geq 1\}}).$$

In order to compute the two terms in the right-hand side of the above equality, we proceed as follows. On $\{J_T=0\}$, there is no jump on $]0,T]$, thus $S_T$ and $\partial_x S_T$ solve some deterministic integral equation. In the examples that we considered in this paper, the solutions of these equations are explicit, so this term is explicitly known. We may use the finite difference method. For the computation of the second term, we use a Monte Carlo algorithm. We simulate a sample $((T_n^k)_{n\in\mathbb{N}},(\Delta_n^k)_{n\in\mathbb{N}})$, $k=1,\dots,M$, of the times and the amplitudes of the jumps and we compute the corresponding $J_t^k$, $S_T^k$ and $H^k_{J_T^k}$. We then write

$$\mathrm{E}(\phi(S_T)H_{J_T}(S_T,\partial_x S_T)\mathbf{1}_{\{J_T\geq 1\}}) \simeq \frac{1}{M}\sum_{k=1}^{M}\phi(S_T^k)H^k_{J_T^k}\mathbf{1}_{\{J_T^k\geq 1\}}.$$

We now compute the Malliavin weights $H^k_{J_T}(S^k_T, \partial_x S^k_T)$ for our examples. One may use Lemma 3.1, but in the particular cases that we discuss here, we have explicit solutions, so direct computations are much easier.

• We first study the diffusion process defined by (4.1). We have the following explicit expression for $S_T$ on $\{J_T = n\}$:

$$S_T = xe^{-rT} + \alpha(1 - e^{-rT}) + \sigma \sum_{j=1}^{n} \Delta_j e^{-r(T-T_j)}. \tag{4.3}$$

We may use integration by parts with respect to the jump amplitudes, to the jump times or to both of them.

∗ Jump amplitudes: $H_{J_T}$ has been calculated in [1]. Since $\Delta_i$ is Gaussian distributed for all $i$, the weight is $\pi(\omega, \Delta_i) = 1$ and on $\{J_T = n\}$,

$$H_n(S_T, \partial_x S_T) = \frac{\sum_{j=1}^{n} e^{rT_j} \Delta_j}{\sigma \sum_{j=1}^{n} e^{2rT_j}}. \tag{4.4}$$

∗ Jump times: suppose that $n \geq 4$ and $J_T = n$. We use the weights $\pi_i(\omega, T_i) = (T_{i+1} - T_i)^\alpha (T_i - T_{i-1})^\alpha$ and we have $\pi'_i = \alpha\delta_{i+1}^{\alpha-1}\delta_i^{\alpha-1}(\delta_{i+1} - \delta_i)$, where $\delta_i = T_i - T_{i-1}$. Then

$$D_i S_T = \sigma\Delta_i r e^{-r(T-T_i)}$$

and

$$L_i S_T = -\sigma r \Delta_i e^{-r(T-T_i)}(r\pi_i + \alpha(\delta_{i+1}\delta_i)^{\alpha-1}(\delta_{i+1} - \delta_i)),$$
$$\sigma_{S_T} = \sum_{i=1}^{n} \pi_i (\sigma r)^2 \Delta_i^2 e^{-2r(T-T_i)}.$$

We define

$$A_j = \alpha(\delta_{j+1}\delta_j)^{\alpha-1}\Delta_j^2 e^{2rT_j},$$
$$B_j = \Delta_j^2 e^{2rT_j}[2r\pi_j + \alpha(\delta_{j+1}\delta_j)^{\alpha-1}(\delta_{j+1} - \delta_j)].$$

Then

$$D_j \sigma_{S_T} = (\sigma r)^2 e^{-2rT}(A_{j-1}\delta_{j-1} - A_{j+1}\delta_{j+2} + B_j).$$

Moreover, $\partial_x S_T = e^{-rT}$ so that $D_i \partial_x S_T = 0$ for all $i = 1, \ldots, n$.

We now have the expression for all of the terms involved in $H_n$ and we obtain

$$\begin{aligned} H_n = {} & \frac{\sum_{i=1}^{n} \Delta_i e^{rT_i}(r\pi_i + \alpha(\delta_{i+1}\delta_i)^{\alpha-1}(\delta_{i+1} - \delta_i))}{\sigma r \hat{\sigma}} \\ & - \frac{\sum_{i=1}^{n} \pi_i \Delta_i e^{rT_i}(A_{i-1}\delta_{i-1} - A_{i+1}\delta_{i+2} + B_i)}{\sigma r \hat{\sigma}^2}, \end{aligned} \tag{4.5}$$

where $\hat{\sigma} = \sum_{i=1}^{n} \pi_i \Delta_i^2 e^{2rT_i}$.

For $n = 1, 2, 3$, we use integration by parts with respect to $\Delta_1$ only. Similar computations then give

$$H_n = \frac{e^{-rT_1}}{\sigma \Delta_1}.$$

• We study the jump diffusion process defined by (4.2). On $\{J_T = n\}$, we have

$$S_T = xe^{rT} \prod_{j=1}^{n} (1 + \sigma \Delta_j).$$

We may not use integration by parts with respect to the jump times because $S_T$ depends on $T_1, \dots, T_k$ by means of $J_t$ only. Therefore, we use the integration by parts formula with respect to the jump amplitudes only. On $\{J_T = n\}$, the Malliavin weight is, in this case (see [1]),

$$H_n(S_T, \partial_x S_T) = \frac{B_\sigma}{\sigma x A_\sigma} + \frac{1}{x} - \frac{2C_\sigma}{x A_\sigma^2}, \tag{4.6}$$

where

$$A_\sigma = \sum_{j=1}^{n} \frac{1}{(1 + \sigma \Delta_j)^2}, \qquad B_\sigma = \sum_{j=1}^{n} \frac{\Delta_j}{(1 + \sigma \Delta_j)}$$

and

$$C_\sigma = \sum_{j=1}^{n} \frac{1}{(1 + \sigma \Delta_j)^4}.$$

**5. Numerical experiments.** In this section, we present several numerical experiments in order to compare the Malliavin approach to the finite difference method. We use the geometrical model and the Vasicek-type model and, in the latter, we use the Malliavin calculus with respect to the amplitudes of the jump and to the jump times. We also look at two types of payoff: call options and digital options.

The comparison is illustrated by some graphs (see Figures 1–3) on one hand and by empirical variance tables (see Tables 1–3) on the other. In Tables 1–3, we give the empirical variance of the two estimators denoted *Var Mall* and *Var Diff* and we compare them. We also include in our tables the value of the volatility $\sigma$ that we use and the corresponding variance of the underlying, denoted by $\mathit{Variance}(S_t)$. We choose the parameter $\sigma$ in the following way.

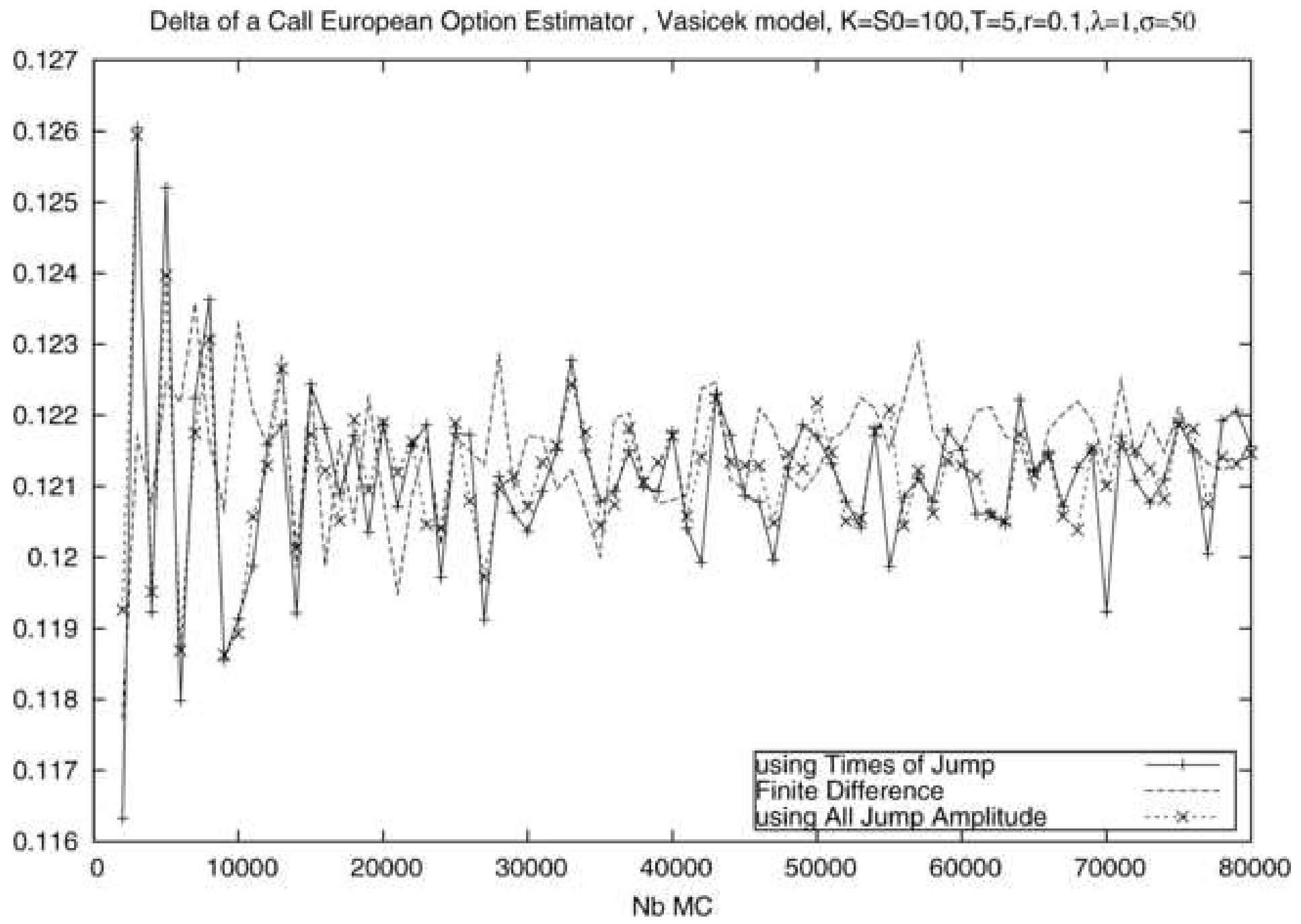

FIG. 1. *Vasicek-type model. Delta of a European call option using Malliavin calculus based on jump times, on jump amplitudes and finite difference method.*

- For the geometrical model, the variance of $S_t$ is

$$Variance(S_t) = x^2 e^{2rt}(e^{\sigma^2 \lambda t} - 1).$$

TABLE 1
*Variance of the Malliavin JT estimator, AJ estimator and of the FD for call option in the Vasicek-type model*

| $\sigma$ | $Variance(S_T)$ | *Var MallJT* | *Var MallAJ* | *Var Diff* |
|---|---|---|---|---|
| 15.8114 | 796.241 | 0.0285123 | 0.0106426 | 0.0300379 |
| 16.6667 | 897.577 | 0.0417219 | 0.0115955 | 0.0298567 |
| 17.6777 | 991.453 | 0.0400695 | 0.013123 | 0.0298904 |
| 18.8982 | 1134.11 | 0.0410136 | 0.0144516 | 0.0299574 |
| 20.4124 | 1313.42 | 0.0433065 | 0.0162378 | 0.029862 |
| 22.3607 | 1584.9 | 0.0400481 | 0.0178726 | 0.0298987 |
| 25 | 1967.53 | 0.0407136 | 0.0202055 | 0.0299007 |
| 28.8675 | 2604.22 | 0.0362728 | 0.0224265 | 0.0299651 |
| 35.3553 | 3961.31 | 0.0343158 | 0.0253757 | 0.0297775 |
| 50 | 7890.4 | 0.0333298 | 0.0287716 | 0.0299749 |

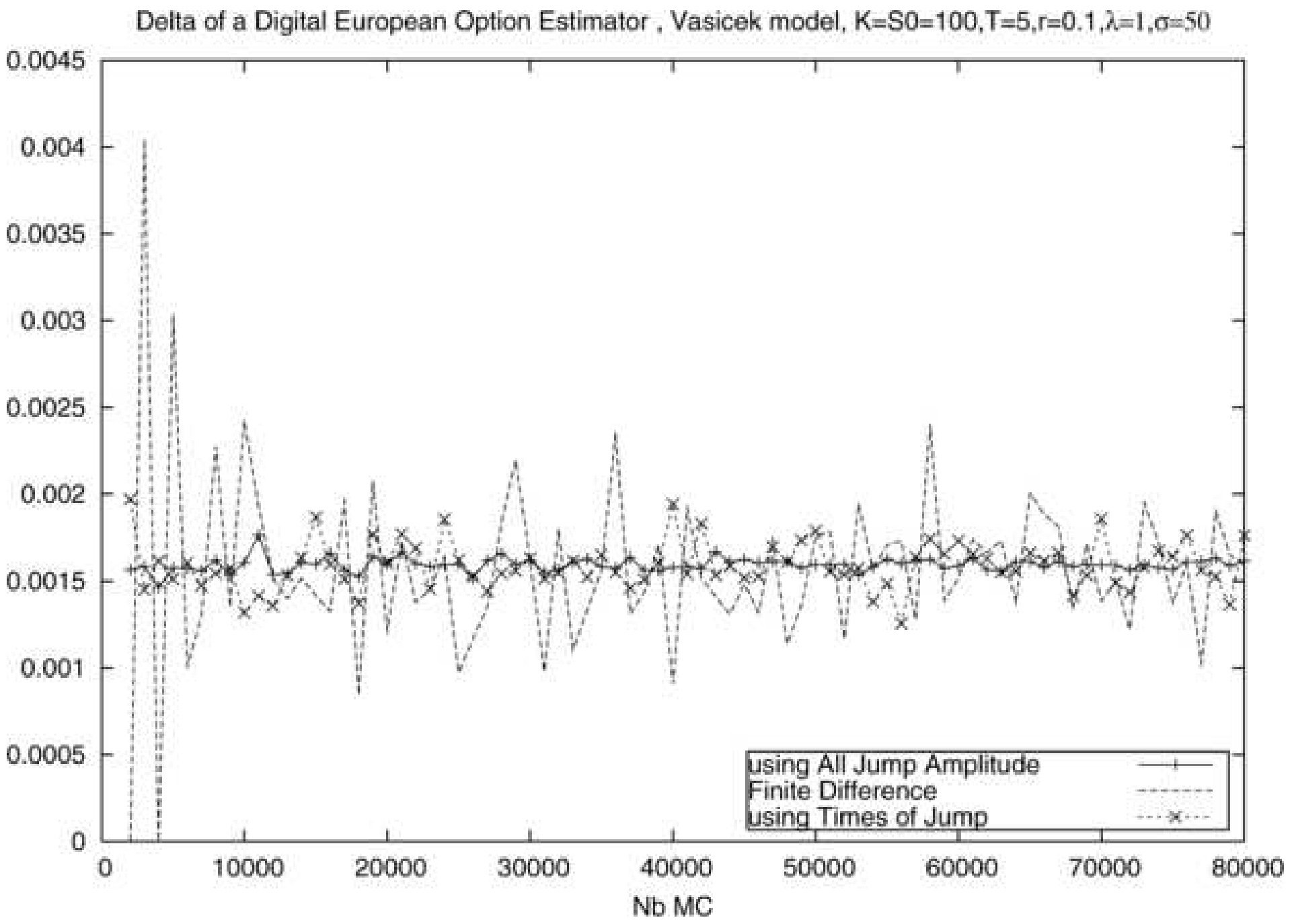

FIG. 2. *Vasicek-type model. Delta of an European digital option using Malliavin calculus based on the jump amplitudes, on the jump times and finite difference method.*

We take $\lambda = 1$, $r = 0.1$, $T = 5$ and $x = 100$. Then for $\sigma \in [0.1, 0.6]$, $1393.69 \leq Variance(S_T) \leq 137264$.

TABLE 2
*Vasicek-type model. Variance of the Malliavin JT estimator, AJ estimator and of the FD for digital option*

| $\sigma$ | $Variance(S_T)$ | Var MallJT | Var MallAJ | Var Diff |
|---|---|---|---|---|
| 15.8114 | 796.241 | 0.00144622 | 7.18878$e$–5 | 0.00514743 |
| 16.6667 | 897.577 | 0.00254652 | 7.3629$e$–5 | 0.00459619 |
| 17.6777 | 991.453 | 0.0018011 | 7.85552$e$–5 | 0.00496369 |
| 18.8982 | 1134.11 | 0.0109864 | 8.14005$e$–5 | 0.00477995 |
| 20.4124 | 1313.42 | 0.00177648 | 8.1627$e$–5 | 0.00386111 |
| 22.3607 | 1584.9 | 0.00152777 | 8.06193$e$–5 | 0.00496369 |
| 25 | 1967.53 | 0.0013786 | 7.94341$e$–5 | 0.0062497 |
| 28.8675 | 2604.22 | 0.00100181 | 7.5835$e$–5 | 0.00551488 |
| 35.3553 | 3961.31 | 0.000617271 | 6.95225$e$–5 | 0.00459619 |
| 50 | 7890.4 | 0.000373802 | 5.64325$e$–5 | 0.00533116 |

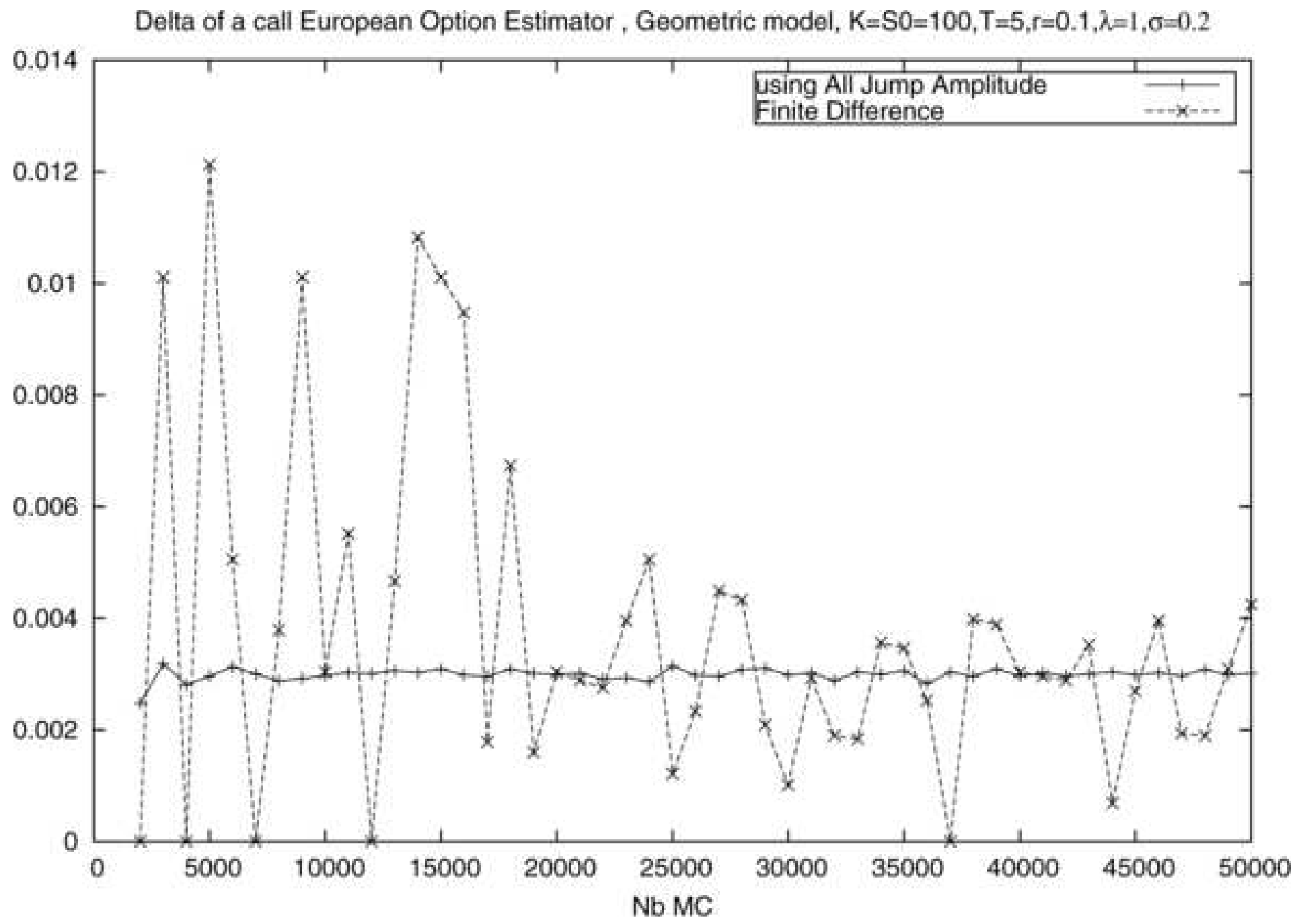

Fig. 3. *Geometrical model. Delta of an European digital option using Malliavin calculus based on the jump amplitudes and finite difference method.*

- For the Vasicek-type model, we have

$$Variance(S_t) = 2\alpha e^{-2rt}(x-\alpha) + \frac{\lambda\sigma^2}{2r}(1-e^{-2rt}).$$

We take $\lambda = 1$, $r = 0.1$, $T = 5$, $\alpha = 10$ and $x = 100$. Then for $\sigma \in [16, 50]$ (note that the practitioners use $\sigma = 20$ to 30 in the short-term rate modeling), $1471.3 \leq Variance(S_T) \leq 8563.69$.

Table 3
*Variance of the Malliavin estimator of the Delta and variance of the FD for digital option*

| $\sigma$ | $Variance(S_T)$ | *Var Mall* | *Var Diff* |
|---|---|---|---|
| 0.1 | 1405.06 | 0.000263425 | 0.00102718 |
| 0.2 | 6183.72 | 0.000917207 | 0.00164801 |
| 0.3 | 16005.5 | 0.000885212 | 0.00117345 |
| 0.4 | 42590.8 | 0.000685313 | 0.0013196 |
| 0.5 | 69018.7 | 0.000531118 | 0.000917399 |
| 0.6 | 130425 | 0.000310461 | 0.0003307 |

In all of our simulations, we have used a variance reduction method based on localization (analogous to the one introduced in [13] and [12]). We use the following abbreviations:

- AJ: amplitude of the jumps.
- JT: jump times.
- FD: finite differences.
- G: geometrical model.
- V: Vasicek-type model.
- Call: call option.
- Dig: digital option.

(G/Dig/AJ) then means that we are in the geometrical model (G) with a digital option (Dig) and we use an algorithm based on the amplitudes of the jumps (AJ). (G/Dig/AJ) versus (G/Dig/FD) means that we compare these two estimators.

5.1. *Numerical results for the Vasicek-type model.* Here, we compare the Delta of European call and digital options obtained by using Malliavin calculus on the one hand and finite difference method on the other hand.

- (V/Call/AJ) versus (V/Call/JT) versus (V/Call/FD).
- (V/Dig/AJ) versus (V/Dig/JT) versus (V/Dig/FD).

In the call options case, both the graph and the variance table show that the Malliavin estimators and the finite difference variances are very close. The methods then give comparable results.

In the digital options case, the algorithm based on Malliavin calculus is significantly better than the one based on finite difference. We also note that the Malliavin estimator using jump amplitudes gives less variance than the one based on jump times.

5.2. *Numerical results for the geometrical model.* These are similar to those of the Vasicek-type model, as we can see for digital options.

(G/Dig/AJ) versus (G/Dig/FD).

5.3. *Conclusions.*

- For a smooth payoff (as the call), the algorithms based on the Malliavin calculus (with respect to the jump times or amplitudes) give comparable results to those based on the finite difference method.
- For a discontinuous payoff (as in the digital options), the algorithms based on Malliavin calculus give significantly better results than those based on the finite difference method.

V. BALLY
UNIVERSITÉ DE MARNE LA VALLÉE
LABORATOIRE D'ANALYSE ET
DE MATHÉMATIQUES APPLIQUÉES
UMR 8050
5 BOULEVARD DESCARTES
CITÉ DESCARTES–CHAMPS SUR MARNE
77454 CEDEX 2
FRANCE
E-MAIL: vlad.bally@univ-mlv.fr

M.-P. BAVOUZET
INRIA ROCQUENCOURT
DOMAINE DE VOLUCEAU-ROCQUENCOURT
PROJET MATHFI
78150 LE CHESNAY
FRANCE
E-MAIL: marie-pierre.bavouzet@inria.fr

M. MESSAOUD
IXIS
47 QUAI D'AUSTERLITZ
75648 PARIS CEDEX 13
AND
INRIA ROCQUENCOURT
PROJET MATHFI
FRANCE
E-MAIL: mmessaoud@ixis-cib.com